\documentclass[11pt,reqno]{amsart}
\usepackage{amsmath,amsthm,amsfonts,amssymb,mathrsfs,bm,graphicx,stmaryrd,hyperref,relsize}
\usepackage{epstopdf}
\usepackage[usenames]{color}
\usepackage[letterpaper,hmargin=1.0in,vmargin=1.0in]{geometry}
\parskip 	\smallskipamount

\newtheorem{theorem}{Theorem}[section]
\newtheorem{lemma}[theorem]{Lemma}
\newtheorem{corollary}[theorem]{Corollary}
\newtheorem{proposition}[theorem]{Proposition}

\def\N{\mathbb{N}}
\def\P{\mathbb{P}}
\def\Z{\mathbb{Z}}
\def\R{\mathbb{R}}

\def\E{\mathbb{E}}

\def\l{\ell}
\def\epsilon{\varepsilon}

\newcommand{\holder}{H\"{o}lder }
\newcommand{\weight}{\mathsf{W}}
\newcommand{\tot}{t_{1,2}}

\newcommand{\midd}{{\rm Mid}}
\newcommand{\south}{{\rm South}}
\newcommand{\north}{{\rm North}}
\newcommand{\west}{{\rm West}}
\newcommand{\strip}{{\rm Strip}}

\newcommand{\high}{{\rm High}}

\begin{document}
\title[Polymer fluctuations and weight profiles in last passage percolation]{Modulus of continuity for polymer fluctuations and weight profiles in Poissonian last passage percolation}

\author{
Alan Hammond}
\address{Departments of Mathematics and Statistics, 
University of California at Berkeley, 
Berkeley, CA, USA.} 
\email{alanmh@stat.berkeley.edu}
\author
{Sourav Sarkar}
\address{Department of Statistics, 
University of California at Berkeley, 
Berkeley, CA, USA.}
\email{souravs@berkeley.edu}

\thanks{The first author is supported by NSF grant DMS-$1512908$. }
\date{}

\date{}

\begin{abstract}
In last passage percolation models, the energy of a path is maximized over all directed paths with given endpoints in a random environment, and the maximizing paths are called {\em geodesics}. The geodesics and their energy can be scaled so that transformed geodesics  cross unit distance and have fluctuations and scaled energy  of unit order. Here we consider Poissonian last passage percolation, a model lying in the KPZ universality class, and refer to scaled geodesics 
as {\em polymers} and 
their scaled energies as {\em weights}. Polymers may be viewed as random functions of the vertical coordinate and, when they are, we show that they have modulus of continuity whose order is at most  $t^{2/3}\big(\log t^{-1}\big)^{1/3}$. 
The power of one-third in the logarithm may be expected to be sharp and in a related problem we show that it is: among polymers in the unit box whose endpoints have vertical separation $t$ (and a horizontal separation of the same order),
the  maximum transversal fluctuation has order $t^{2/3}\big(\log t^{-1}\big)^{1/3}$.
Regarding the orthogonal direction, in which growth occurs, we show that, when one endpoint of the polymer is fixed at $(0,0)$ and the other is varied vertically  over $(0,z)$, $z\in [1,2]$, the  resulting random weight profile has sharp modulus of continuity of order $t^{1/3}\big(\log t^{-1}\big)^{2/3}$. 
In this way, we identify exponent pairs of $(2/3,1/3)$
and $(1/3,2/3)$ in power law and polylogarithmic correction, respectively for polymer fluctuation, and polymer weight under vertical endpoint perturbation.
The two exponent pairs describe~\cite{H121,H122,H123} the fluctuation of the boundary separating two phases in subcritical planar random cluster models.
\end{abstract}

\maketitle
\vspace{-.1in}
         
\tableofcontents

\section{Introduction} 
In 1986, Kardar, Parisi, and Zhang \cite{KPZ86} predicted universal scaling behaviour for many planar random growth processes, including first and last passage percolation as well as corner growth processes, though rigorous validation has been subsequently provided for only a handful of them. In such models, fluctuation in the direction of growth is governed by an exponent of one-third, with this fluctuation enduring on a scale governed by an exponent of two-thirds in the orthogonal, or transversal, direction.

 Poissonian last passage percolation illustrates these effects. We will define it shortly, since it is our object of study; briefly, the model specifies a growth process whose height at a given moment is the maximum number of points (or the {\em energy}) obtainable in a directed path  through a planar Poisson point process. Baik, Deift and Johansson \cite{BDJ99} established the $n^{1/3}$-order fluctuation of  the maximum number of Poisson points on an increasing path from $(0,0)$ to $(n,n)$, deriving the GUE Tracy-Widom distributional limit of the scaled energy. Later Johansson \cite{J00} proved the transversal fluctuation exponent of two-thirds in this model.
These are exactly solvable models, for which certain exact distributional formulas are available, and the derivations of these formulas typically employ deep machinery from algebraic combinatorics or random matrix theory. 
It is interesting to study geometric properties of universal KPZ objects by approaches that, while they are reliant on certain integrable inputs, are probabilistic in flavour: for example, \cite{BSS14},\cite{BSS17++} and \cite{BSS17+} are recent results and applications concerning geometric properties of last passage percolation paths. 

It is rigorously understood, then, that last passage percolation paths experience fluctuation in their energy and transversal fluctuation governed by scaling exponents of one-third and two-thirds. It is very natural to view such paths via the lens of scaled coordinates, in which transversal fluctuation and path energy each has unit order. We will be more precise very shortly, when suitable notation has been introduced, but for now we mention that our aim in this article is to refine rigorous understanding of the magnitude and geometry of fluctuation in last passage percolation paths. We shall call the scaled geodesics {\em polymers}, and refer to the scaled energy as {\em weight}. We will see that polylogarithmic corrections to the scaled laws implied by the exponents of one-third and two-thirds arise when we consider natural geometric problems concerning the weights and the maximum fluctuation among polymers in a unit order region. 
The techniques for verifying our claims will employ geometric and probabilistic tools rather than principally integrable ones, since problems involving maxima as both endpoints of a last passage percolation path are varied are not usually amenable to integrable techniques.

\subsection{Model definition and main results}\label{ss:defres}
Let $\Pi$ be a homogeneous rate one Poisson point process (PPP) on $\R^2$. We introduce a partial order on $\R^2$:  $(x_1,y_1)\preceq (x_2,y_2)$ if and only if $x_1\leq x_2$ and $y_1\leq y_2$.  For $u\preceq v$, $u,v\in \R^2$, an increasing path $\gamma$ from $u$ to $v$ is a piecewise affine path, viewed as a subset of $\R^2$, 
that joins points $u=u_0\preceq u_1\preceq u_2 \preceq \ldots \preceq u_k = v$ such that $u_i\in \Pi$ for $i\in \llbracket 1,k-1\rrbracket$. Here and later, $\llbracket a,b\rrbracket$ for $a,b\in \Z$ with $a\leq b$ denotes the integer interval $\{a,\cdots,b\}$. Also let $|\gamma|$ denote the {\em energy} of $\gamma$, namely the number of points in $\Pi\setminus\{v\}$ that lie on $\gamma$; (the last vertex is excluded from the definition of energy so that the sum of the energies of two paths equals the energy of the concatenated path, as we will see in Section \ref{ss:pocat}). Then we define the last passage time from $u$ to $v$, denoted by $X_u^v$, to be the maximum of $|\gamma|$ as $\gamma$ varies over all increasing paths from $u$ to $v$. Any such maximizing path is called a geodesic. There may be several such, but if $\Gamma_u^v$ denotes any one of them, we have
\begin{equation}\label{e:defE}
X_u^v=|\Gamma_u^v| \, .
\end{equation}
Note that, in this notation, the
starting and ending points of the geodesic, $u$ and $v$, are assigned subscript and superscript placements. 
We will often use this convention, including in the case of the scaled coordinates that we will introduce momentarily.

When $u\preceq v$, any geodesic from $u$ to $v$
may be viewed as a function of its horizontal coordinate, since it contains a vertical line segment with probability zero. The operations of maximum and minimum may be applied to any pair of such geodesics, and the results are also geodesics.
For this reason, we may speak unambiguously of $\Gamma_{u}^{\leftarrow;v}$, the uppermost geodesic between $u$ and $v$, and of $\Gamma^{\rightarrow;v}_{u}$, the lowermost geodesic between $u$ and $v$. 
(The notation $\leftarrow$ and $\rightarrow$ is compatible with these two paths being  equally well described as the leftmost and rightmost geodesics. This choice of notation also anticipates the form of these paths when viewed in the scaled coordinates that we are about to introduce.)
When the endpoints are $(0,0)$ and $(n,n)$, we will call these geodesics $\Gamma_n^{\leftarrow}$ and~$\Gamma_n^{\rightarrow}$.

\subsubsection{Introducing scaled coordinates}
 We rotate the plane about the origin counterclockwise by  $45$ degrees, squeeze the vertical coordinate  by a factor $2^{1/2}n$  and the horizontal one by $2^{1/2}n^{2/3}$, thus setting  
 \begin{equation}\label{e:deftrans}
 T_n:(x,y)\mapsto\left(2^{-1}n^{-2/3}(x-y),2^{-1}n^{-1}(x+y)\right) .
 \end{equation}

The horizontal line at vertical coordinate $t$
is the image under $T_n$
of the anti-diagonal line through $(nt,nt)$.  It is easy to see that, for $(x,t)\in \R^2$, $T_n^{-1}(x,t)=(nt+xn^{2/3},nt-xn^{2/3})$.
 
\begin{figure}
\centering

\includegraphics[scale=.75]{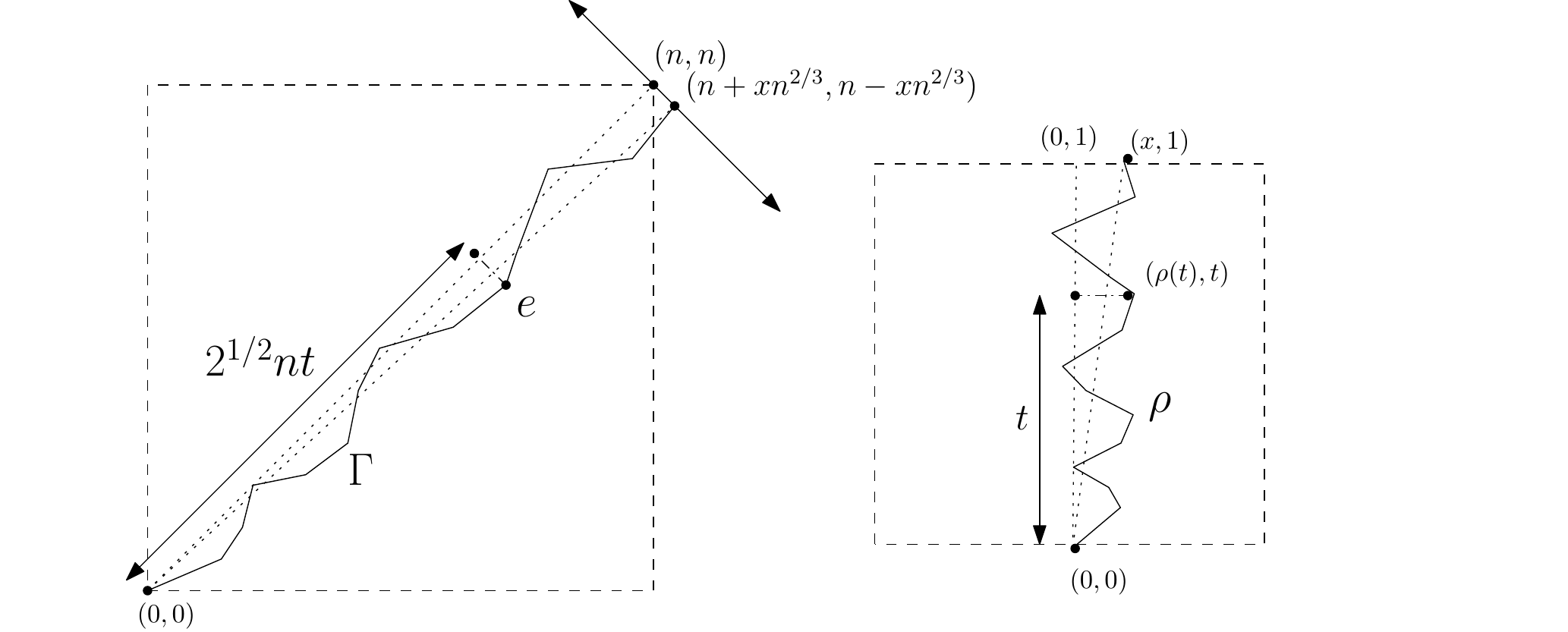} 
\caption{The scaling map $T_n$ applied to the left figure produces the figure on the right. The point $e$ in the geodesic $\Gamma$ is the preimage of the point $(\rho(t),t)$ in the polymer $\rho$. }
\label{f:scaling}
\end{figure}

Paths that are the image of geodesics under $T_n$ will be called {\em polymers}; we might say $n$-polymers, but the suppressed parameter will always be $n$. 
Geodesics from $(0,0)$ to $(n,n)$ transform to polymers $(0,0)$ to $(0,1)$. 
Figure~\ref{f:scaling}
depicts a geodesic~$\Gamma$ and its image polymer~$\rho$. 
The polymer between planar points~$u$ and~$v$
that is the image of the uppermost geodesic given the preimage endpoints will be denoted by
$\rho^{\leftarrow;v}_{n;u}$, and, naturally enough, called the leftmost polymer from $u$ to $v$.
The rightmost polymer from $u$ to $v$ is the image of the corresponding lowermost geodesic and will be denoted by $\rho^{\rightarrow;v}_{n;u}$. The simpler notation $\rho^{\leftarrow}_n$ and $\rho^{\rightarrow}_n$ will be adopted when $u = (0,0)$ and $v = (0,1)$. When $u=(x_1,t_1), v=(x_2,t_2)$, with $x_1,x_2,t_1,t_2\in \R$, $t_1<t_2$, such that $T_n^{-1}(x_1,t_1)\preceq T_n^{-1}(x_2,t_2)$, 
we will, when it is convenient, regard any polymer $\rho$ from $u$ to $v$ as a function of its vertical coordinate: that is, for $t \in [t_1,t_2]$, $\rho(t)$ will denote the unique point such that $(\rho(t),t)\in \rho$.
(This definition makes sense since an increasing  path can intersect any anti-diagonal at most once.) We regard the vertical coordinate as time, as the $t$-notation suggests, and will sometimes refer to the interval $[t_1,t_2]$ as the {\em lifetime} of the polymer. In particular, when $t_1=0,t_2=1$,
writing $C[0,1]$ for the space of
continuous real-valued functions on $[0,1]$ (equipped for later purposes with the topology of uniform convergence), we may thus view $\rho=\{\rho(t)\}_{t \in [0,1]}$ as an element of $C[0,1]$.

\subsubsection{Condition for existence of polymers}\label{sss:comp} For $u=(x_1,t_1),v=(x_2,t_2)$ with $x_1,x_2,t_1,t_2\in \R$, $t_1<t_2$, we have that $T_n^{-1}(u)=(nt_1+x_1n^{2/3},nt_1-x_1n^{2/3})$ and $T_n^{-1}(v)=(nt_2+x_2n^{2/3},nt_2-x_2n^{2/3})$. Thus $T_n^{-1}(u)\preceq T_n^{-1}(v)$ is and only if $|x_1-x_2|<n^{1/3}(t_2-t_1)$. Indeed, we will write $u\overset{n}{\preceq}v$ to mean that $|x_1-x_2|<n^{1/3}(t_2-t_1)$; this condition ensures that polymers exist between the endpoints $u$ and $v$.

The first of our three main results shows that polymers, viewed as functions of the vertical coordinate, enjoy modulus of continuity of order $t^{2/3} \big( \log t^{-1} \big)^{1/3}$.

\begin{theorem}\label{t:main1}
\begin{itemize}
\item[(a)] The sequence $\{\rho_n^{\leftarrow}\}_{n\in \N}$ is tight in $(C[0,1],\|\cdot\|_\infty)$. 
\item[(b)]There exists  a constant $C > 0$ such that, for the weak limit $\rho_*^{\leftarrow}$ of any weakly converging subsequence of $\{\rho_n^{\leftarrow}\}_{n\in \N}$, almost surely,
\begin{equation}\label{e:modcon}
\limsup_{t\searrow 0}\sup_{0\leq z\leq 1-t}{t^{-2/3}\big(\log t^{-1}\big)^{-1/3}}{|\rho_*^{\leftarrow}(z+t)-\rho_*^{\leftarrow}(z)|}\leq C \, .
\end{equation}
The same result holds for the rightmost polymer.
\end{itemize}
\end{theorem}
Note that  the constant $C$ does not depend on the choice of the weakly converging subsequence.


The exponent pair $(2/3,1/3)$ for power law and polylogarithmic correction is thus demonstrated to hold in an upper bound on polymer fluctuation. 
We believe that a lower bound holds as well, in the sense that the limit infimum counterpart to 
  \eqref{e:modcon} is positive.  
A polymer is an object specified by a global constraint, and it by no means clearly enjoys independence properties as it traverses disjoint regions, even though the underlying Poisson randomness does. In order to demonstrate the polymer fluctuation lower bound, this subtlety would have to be addressed.
We choose instead to demonstrate that the exponent pair $(2/3,1/3)$ describes polymer fluctuation by  proving a lower bound of this form for the maximum fluctuation witnessed among a natural class of short polymers in a unit region. This alternative formulation offers a greater supply of independent randomness.

Indeed, we now specify a notion of {\em maximum transversal fluctuation} over a collection of short polymers. Fix any two points $u=(x_1,t_1),v=(x_2,t_2)$ such that $t_2>t_1$. 
Let $\Phi_{n;u}^v$ denote the set of all polymers $\rho$ from $u$ to $v$. Let $\l_{u}^v$ denote the planar line segment that joins $u$ and $v$; extending an abuse of notation that we have already made, we write  $\l_u^v(t)$ for the unique point such that $(\l_u^v(t),t)\in \l_u^v$, where $t\in [t_1,t_2]$. Then, for any polymer $\rho$, the transversal fluctuation 
$\mathrm{TF}(\rho)$ of $\rho$ is specified to be
\begin{equation}\label{e:TFfixedpol}
\mathrm{TF}(\rho):=\sup_{t\in [t_1,t_2]}|\rho(t)-\l_u^v(t)|,
\end{equation}
and the transversal fluctuation between the points $u$ and $v$ to be
\begin{equation}\label{e:TFbetwnpoints}
\mathrm{TF}_{n;u}^v:=\max_{\rho\in \Phi_{n;u}^v}\mathrm{TF}(\rho)=\max\left\{\mathrm{TF}(\rho_{n;u}^{\leftarrow;v}),\mathrm{TF}(\rho_{n;u}^{\rightarrow;v})\right\}.
\end{equation}

Also, let
$$
\mathsf{InvSlope}_{(x_1,t_1)}^{(x_2,t_2)}=\frac{x_2-x_1}{t_2-t_1} 
$$
denote the reciprocal of the slope of the interpolating line. Since $t_2>t_1$, $\mathsf{InvSlope}_{(x_1,t_1)}^{(x_2,t_2)} \in \R$.

Now fix some large constant $\psi>0$. Then, for any fixed parameter $t\in (0,1]$ and any $n\in \N, n>\psi^3$, we define the set of {\em admissible endpoint pairs}  
\begin{eqnarray}\label{e:defComp}
\mathrm{AdEndPair}_{n,\psi}(t):=\Big\{((x_1,t_1),(x_2,t_2)):t_2-t_1\in (0,t] \, , \, \left|\mathsf{InvSlope}_{(x_1,t_1)}^{(x_2,t_2)}\right|\leq \psi,\nonumber \\
 x_1,x_2\in [-1,1] \, , \,  t_1,t_2\in [0,1]\Big\}.
\end{eqnarray}
Since $n>\psi^3$,
\[|x_2-x_1|n^{2/3}\leq \psi(t_2-t_1)n^{2/3}<(t_2-t_1)n\,.\]
Recalling the notation at the start of Subsection~\ref{sss:comp}, we thus have $(x_1,t_1)\overset{n}\preceq (x_2,t_2)$, so that polymers do exist between such endpoint pairs.

We then define 
\begin{equation}\label{e:defMTF}
\mathrm{MTF}_n(t) = \mathrm{MTF}_{n,\psi}(t):=\sup\left\{\mathrm{TF}_{n;u}^{v}:(u,v)\in \mathrm{AdEndPair}_{n,\psi}(t) \right\},
\end{equation}
so that $\mathrm{MTF}_n(t)$ is the {\em maximum transversal fluctuation} over polymers between all endpoint pairs at vertical distance at most $t$ such that the slope of the interpolating line segment is bounded away from being horizontal; (we suppress the parameter $\psi$ in the notation). Our second theorem demonstrates that the exponent pair $(2/3,1/3)$ governs this maximum traversal fluctuation.


\begin{theorem}\label{t:flucshortpol} There exist $\psi$-determined constants $0<c<C<\infty$  such that 
\[\liminf_n \, \P\left( \, t^{-2/3}\big(\log t^{-1}\big)^{-1/3} \mathrm{MTF}_n(t)  \in [c,C] \, \right)  \to 1 \quad \mbox{as } \, \, {t\searrow 0}\,.\]

\end{theorem}

\subsubsection{Scaled energies are called weights}
It is natural to scale the energy of a geodesic when we view the geodesic as a polymer after scaling. Scaled energy will be called {\em weight}
and specified so that it is of unit order for polymers that cross unit-order distances. For $t_1<t_2$, let $\tot$ denote $t_2-t_1$; (this is a notation that we will often use). Let $(x,t_1),(y,t_2)\in \R^2$ be such that $|x-y|<\tot n^{1/3}$. (This condition ensures that $(x,t_1)\overset{n}\preceq (y,t_2)$, so that polymers exist between this pair of points.) Since $T_n^{-1}((x,t_1))=( nt_1 + xn^{2/3},nt_1 - xn^{2/3})$ and $T_n^{-1}((y,t_2))=( nt_2 + yn^{2/3},nt_2 - yn^{2/3})$, it is natural to define the scaled energies, which we call \textit{weights}, in the following way. 
Define
\begin{equation}\label{e:defweight}
\weight_{n;(x,t_1)}^{(y,t_1)} = n^{-1/3} \bigg( X_{( nt_1 + n^{2/3}x, nt_1 - n^{2/3}x)}^{(nt_2 + n^{2/3}y, nt_2 - n^{2/3}y)} \Big) \, - \, 2 n \tot \bigg) \, .
\end{equation}

Because of translation invariance of the underlying Poisson point process, $\tot$ is a far more relevant parameter than $t_1$ or $t_2$. The notation on the left-hand side of~(\ref{e:defweight})
is characteristic of our presentation in this article: a scaled object is being denoted, with planar points $(\cdot,\cdot)$ in the subscript and superscript indicating starting and ending points.

\subsubsection{A continuous modification of the weight function}
For the statement of our third theorem, we prefer to make an adjustment to the polymer weight to cope with a minor problem
concerning discontinuity of geodesic energy under endpoint perturbation. 
For $n\in \N$, define $X_n:[1,2]\mapsto [0,\infty)$,
\[X_n(t):=X_{(0,0)}^{(nt,nt)} \, . \]
Observe that $X_n(t)$ is integer-valued, non-decreasing, right continuous and has almost surely a finite number of jump discontinuities. 
Let $d_0 = 1$ and $d_m = 2$.
Record in increasing order the points of discontinuity of $X_n$ 
as a list   $\big( d_1,d_2,\cdots,d_{m-1} \big)$. We specify a {\em modified} and continuous form of the function $X_n$ by linearly interpolating it between these points of discontinuity, setting 
\[X_n^{\mbox{\tiny{mod}}}(t):=X_n(d_i)+(t-d_i)(d_{i+1}-d_i)^{-1}\big(X_n(d_{i+1})-X_n(d_i)\big), \mbox{ for } t \in [d_i,d_{i+1}], \]
for $i=1,2,\cdots,m-1$. Because almost surely  no two points in a planar Poisson point process share either their horizontal or vertical coordinate, $X_n(d_{i+1})-X_n(d_i)=1$ for all $i$. Thus, for all $t\in [1,2]$,
\begin{equation}\label{e:Xmod}
X_n(t)\leq X_n^{\mbox{\tiny{mod}}}(t)\leq X_n(t)+1\,.
\end{equation}
 Now define the \textit{modified weight} function $\mathsf{Wgt}_n:[1,2]\mapsto \R$ for polymers from $(0,1)$ to $(\cdot,1)$:
\begin{equation}\label{e:defWgt}
\mathsf{Wgt}_n(t):=n^{-1/3}\left(X_n^{\mbox{\tiny{mod}}}(t)-2nt\right)
\end{equation}
Because of \eqref{e:Xmod},
\begin{equation}\label{e:Wmod}
\left| \,  \mathsf{Wgt}_n(t)-\weight_{n;(0,0)}^{(0,t)} \,  \right| \leq n^{-1/3} \, .
\end{equation}


By construction, $\mathsf{Wgt}_n$ sending $t \in [1,2]$  to  $\mathsf{Wgt}_n(t)$ is an element of $C[1,2]$, the space of continuous functions on $[1,2]$; (similarly to before, this space will be  equipped  with the topology of uniform convergence).

Our third main result demonstrates that the exponent pair $(1/3,2/3)$ offers a description of the modulus of continuity of polymer weight when one endpoint is varied vertically. 
\begin{theorem}\label{t:holdweight} The sequence $\{\mathsf{Wgt}_n\}_{n\in \N}$ is tight in $(C[1,2],\|\cdot\|_\infty)$. There exist constants $0<c<C<\infty$ such that, for the weak limit $\mathsf{Wgt}_*$ of any weakly converging subsequence of $\{\mathsf{Wgt}_n\}_{n\in \N}$, almost surely
\begin{eqnarray}\label{e:wgtholdbounds}
c&\leq &\liminf_{t\searrow 0}\sup_{1\leq z\leq 2-t} \,  t^{-1/3}\big(\log t^{-1}\big)^{-2/3} \, \Big\vert \mathsf{Wgt}_*(z+t)-\mathsf{Wgt}_*(z) \Big\vert \\
&\leq & \limsup_{t\searrow 0}\sup_{1\leq z\leq 2-t}    \, t^{-1/3}\big(\log t^{-1}\big)^{-2/3}  \, \Big\vert \mathsf{Wgt}_*(z+t)-\mathsf{Wgt}_*(z) \Big\vert \, \leq \, C \, .\nonumber
\end{eqnarray}

\end{theorem}

Note that, as in Theorem~\ref{t:main1}, the constants $c$ and $C$ do not depend on the choice of weak limit point or converging subsequence.

Beyond these three theorems, 
we present a proposition, which is needed for the proof of Theorem~\ref{t:flucshortpol} but which also has independent interest.  The maximum  fluctuation of any geodesic joining $(0,0)$ and $(n,n)$ around the interpolating line has probability at most $e^{-c k^3}$ of exceeding $k n^{2/3}$. This upper bound has  essentially been obtained in \cite[Theorem $11.1$ and Corollary $11.7$]{BSS14}, though we will state and prove this result, with the power of three in the exponent inside the exponential, as Theorem \ref{t:transversal}. Our next proposition is the  matching lower bound, stated using scaled coordinates. Observe from \eqref{e:TFfixedpol} that, for any polymer $\rho$ between $(0,0)$ and $(0,1)$, $\mathrm{TF}(\rho)=\sup_{y\in [0,1]}|\rho(y)|$. Also recall that $\Phi_{n;(0,0)}^{(0,1)}$ is the set of all polymers from $(0,0)$ to $(0,1)$. 

\begin{proposition}\label{p:lowbndTF} There exist positive constants $c^*$, $n_0$, $s_0$ and $\alpha_0$ such that, for all $t_1,t_2$ with $\tot=t_2-t_1>0$ and all $n\tot\geq n_0$ and $s\in [s_0,\alpha_0(n\tot)^{1/3}]$, 
$$
\P\left( \min \left\{ \mathrm{TF}(\rho) : \rho\in \Phi_{n;(0,t_1)}^{(0,t_2)} \right\}  \geq  s \tot^{2/3} \right) \,  \geq \, \exp \big\{-c^*s^3 \big\} \, .
$$
\end{proposition}
\subsection{A few words about the proofs} The main ingredients in the proofs of Theorem~\ref{t:main1} and Theorem~\ref{t:flucshortpol} are the estimates from integrable probability assembled in Section~ \ref{s:estinteg} and a {\em polymer ordering} property elaborated in Lemma~\ref{l:porder} that propagates control on polymer fluctuation among polymers whose endpoints lie in a discrete mesh to all polymers in the region of this mesh. The basic tools in the proof of the upper bound in Theorem \ref{t:holdweight} and that of Proposition~\ref{p:lowbndTF} are surgical techniques and comparisons of the weights of polymers, and are reminiscent of the techniques developed and extensively used in \cite{BSS14} and \cite{BSS17++}.

\subsection{Phase separation and KPZ}

Certain random models manifest the scaling exponents of KPZ universality and some of its qualitative features, without exhibiting the richness of behaviour of models in this class. For example, the least convex majorant of the stochastic process $\R \to \R: x \to B(x) - t^{-1} x^2$ is comprised of planar line segments, or facets, the largest of which in a compact region has length of order $t^{2/3 + o(1)}$ when $t > 0$ is high; and the typical deviation of the process from its majorant scales as $t^{1/3 + o(1)}$.

Some such models form a testing ground for KPZ conjectures. Phase separation concerns the form of the boundary of a droplet of one substance suspended in another. When supercritical bond percolation on $\Z^2$ is conditioned on the cluster (or droplet) containing the origin being finite and large, namely of finite size at least $n^2$, with $n$ high, the interface at the boundary of this cluster is expected to exhibit KPZ scaling characteristics, with the scaling parameter $n$ playing a comparable role to $t$ in the preceding example. Indeed, in~\cite{H121,H122,H123}, a surrogate of this interface, expressed in terms of the random cluster model, was investigated. 
The maximum length of the facets that comprise the boundary of the interface's convex hull was proved to typically have the order $n^{2/3} \big( \log n \big)^{1/3}$, while the maximum local roughness, namely the maximum distance from a point on the interface to the convex hull boundary, was shown to be of the order of $n^{1/3} \big( \log n \big)^{2/3}$.

Viewed in this light, the present article validates for the KPZ universality class the implied predictions: 
that exponent pairs of $(1/3,2,3)$ and $(2/3,1/3)$
for power-law and logarthmic-power govern maximal polymer weight change under vertical endpoint displacement and maximal transversal polymer fluctuation. 

In a natural sense, these two exponent pairs are accompanied by a third, namely $(1/2,1/2)$,
for interface regularity. In the example of parabolically curved Brownian motion, $x \to B(x) - x^2 t^{-1}$, the modulus of continuity of the process on $[-1,1]$ is easily seen to have the form $s^{1/2} \big( \log s^{-1} \big)^{1/2}$, up to a random constant, and uniformly in $t \geq 1$. In KPZ, this assertion finds a counterpart when it is made for the Airy$_2$ process, which offers a limiting description in scaled coordinates of the weight of polymers of given lifetime with first endpoint fixed. This assertion has been proved in~\cite[Theorem~$1.11(1)$]{H16}. 
 Recently, for a very broad class of initial data, the polymer weight profile was shown in~\cite[Theorem~ $1.2$]{H17}  to have a modulus of continuity of the order of $s^{1/2}\big( \log s^{-1}\big)^{2/3}$, uniformly in the scaling parameter and the initial condition.


\subsection{Organization} 
We continue with two sections that offer basic general 
tools. The first,
Section~\ref{s:estinteg}, provides useful estimates available from the integrable probability literature. Then, in Section \ref{s:polorder}, we state and prove the polymer ordering lemmas and some other basic results, which are essential tools in the proofs of the main theorems. 

The remaining four sections, \ref{s:2/3,1/3} -- \ref{s:last}, contain the main proofs. Consecutively, these sections are devoted to proving:
\begin{itemize}
\item the polymer \holder continuity upper bound Theorem \ref{t:main1};
\item  the modulus of continuity for maximum transversal fluctuation over short polymers, Theorem \ref{t:flucshortpol}, subject to assuming Proposition \ref{p:lowbndTF};
\item \holder continuity for the polymer weight profile, Theorem \ref{t:holdweight};
\item and the lower bound on transversal polymer fluctuation, Proposition~\ref{p:lowbndTF}.
\end{itemize} 

We will stick to scaled coordinates in the 
results'
statements and, except in Section~\ref{s:estinteg}, in their proofs. A bridge between scaled coordinates and the original ones is offered in this next section, in whose proofs we use the scaling map $T_n$ from \eqref{e:deftrans} and weight function $\weight$ from \eqref{e:defweight} to transfer   unscaled results to their scaled counterparts.

\section{Scalings and estimates from integrable probability}\label{s:estinteg}

In this section, we assemble some results from integrable probability. Most of these results were derived in terms of unscaled coordinates in \cite{BSS14} and \cite{BSS17++}. Point-to-point estimates of last passage percolation geodesics were used crucially in \cite{BSS14} to resolve the ``slow-bond" conjecture, and in \cite{BSS17++} to show the coalescence of nearby geodesics, and those estimates will be crucially employed in this paper as well. We state the results in scaled coordinates, and the proofs detail how to obtain these statements from their unscaled versions available in the literature. In going from the unscaled to scaled coordinates, we shall use the definitions of the scaling map in \eqref{e:deftrans} and the weight in \eqref{e:defweight}. First we observe some simple relations between the different scaled versions of these quantities that will be used in the proofs of the theorems in this section.
\medskip
\\
\noindent{\bf The scaling principle.} Because of translation invariance and the definition \eqref{e:deftrans}, it is easy to see that for any $x,y,t_1,t_2\in \R$ with $\tot=t_2-t_1>0$ and $(x,t_1)\overset{n}{\preceq}(y,t_2)$ (see Subsection \ref{sss:comp}), for any $\theta\in [0,1]$,
\begin{equation}\label{e:screl}
\rho_{n;(x,t_1)}^{\leftarrow;(y,t_2)}(t_1+\theta \tot)\overset{d}{=}\tot^{2/3}\rho_{n;(x\tot^{-2/3},0)}^{\leftarrow;(y\tot^{-2/3},1)}(\theta)\overset{d}{=}\tot^{2/3}\rho_{n;(0,0)}^{\leftarrow;((y-x)\tot^{-2/3},1)}(\theta)\,.
\end{equation}
The same statement holds for the rightmost polymers as well. Here and throughout $\overset{d}{=}$ denotes that the two random variables on either side have the same distribution.
We will sometimes call the displayed assertion the {\em scaling principle}. 

Also by translation invariance and the definition of weight in \eqref{e:defweight}, it follows that
\begin{equation}\label{e:eqdist}
\tot^{-1/3}\weight_{n;(x,t_1)}^{(y,t_2)}\overset{d}{=}\weight_{n\tot;(xt_{1,2}^{-2/3},0)}^{(yt_{1,2}^{-2/3},1)}\overset{d}{=}\weight_{n\tot;(0,0)}^{((y-x)t_{1,2}^{-2/3},1)}\,.
\end{equation}
\medskip
\\
\noindent{\bf Boldface notation for applying results.}
In our proofs, we will naturally often be applying tools such as those stated in this section. Sometimes the notation of the tool and of the context of the application will be in conflict. To alleviate this conflict, we will use boldface notation when we specify the values of the parameters of a given tool in terms of quantities in the context of the application. We will first use this notational device shortly, in one of the upcoming proofs.

\medskip


The next theorem was proved in \cite{BDJ99}.
\begin{theorem}
\label{t:BDJ}
As $n\rightarrow \infty$,
\[\weight_{n;(0,0)}^{(0,1)}\Rightarrow F_{TW},\]
where the convergence is in distribution and $F_{TW}$ denotes the GUE Tracy-Widom distribution. 
\end{theorem}
For a definition of the GUE Tracy-Widom distribution, also called the $F_2$ distribution, see \cite{BDJ99}. 

Moderate deviation inequalities for this centred and scaled polymer weight will be important. Such inequalities follow immediately from \cite[Theorem ~$1.3$]{LM01}, \cite[Theorem~$1.2$]{LMS02} and \eqref{e:eqdist}. These are  essential inequalities, used repeatedly in this paper. In fact, it should be possible to recover the theorems of this paper for other integrable models for which such moderate deviation estimates are known.
 
 \begin{theorem}\label{t:moddev} There exist positive constants $c,s_0$ and $n_0$ such that, for all $t_1<t_2$ with $n\tot >n_0$ and $s>s_0$, 
 \[\P\left(\tot^{-1/3}\weight_{n;(0,t_1)}^{(0,t_2)}\geq s\right)\leq e^{-cs^{3/2}},\]
 and
 \[\P\left(\tot^{-1/3}\weight_{n;(0,t_1)}^{(0,t_2)}\leq -s\right)\leq e^{-cs^{3/2}}.\]
\end{theorem}  

Also, we
shall need not just tail bounds for weights of point to point polymers, but uniform tail bounds on polymer weights whose endpoints vary over fixed unit order intervals.
The unscaled version of this theorem follows from \cite[Propositions~$10.1$ and~$10.5$]{BSS14}.  

\begin{theorem}\label{t:unifmoddev}
There exist $C,c \in (0,\infty), C_0\in (1,\infty)$ and $n_0 \in \N$
such that, for all $t_1<t_2$ with $n\tot \geq n_0$,   $s\in [0,10(n\tot)^{2/3}]$, $A=C_0^{-1}s^{1/4}n^{1/6}\tot^{5/6}$ and $I$ and $J$ intervals of length at most $\tot^{2/3}$ that are contained  in $\left[-A,A\right]$, 
$$
 \P \bigg( \, \sup_{x \in I, y \in J} \Big\vert \tot^{-1/3} \weight_{n;(x,t_1)}^{(y,t_2)} +  \tot^{-4/3} (x-y)^2 \Big\vert > s \, \bigg) \, \leq \,  C \exp \big\{ -c s^{3/2} \big\} \, .
$$
\end{theorem} 

\noindent{\bf Proof.} First we prove the theorem when $t_1=0$ and $t_2=1$ by invoking the unscaled version of this theorem from \cite{BSS14}. At the end we prove Theorem \ref{t:unifmoddev} for general $t_1<t_2$. Observe that $|x-y|<2C_0^{-1}s^{1/4}n^{1/6}\tot^{5/6}<2^{-1}n\tot^{1/3}$ for $C_0>2\cdot10^{1/4}$ since $s\leq 10(n\tot)^{2/3}$. This ensures that $\weight_{n;(x,t_1)}^{(y,t_2)}$ is well defined.

Let $u=T_n^{-1}(x,0)=(xn^{2/3},-xn^{2/3})$ and $v=T_n^{-1}(y,1)=(n+yn^{2/3},n-yn^{2/3})$. If $S_{u,v}$ denotes the slope of the line segment joining $u$ and $v$, then $|x-y|<2^{-1}n$ ensures that $3^{-1}<S_{u,v}<3$. Then, using the first order estimates (see \cite[Corollary $9.1$]{BSS14}) and a simple binomial expansion giving $|(1-x)^{1/2}-(1-2^{-1}x)|\leq C_1x^2$ for $x\in (-1,1)$, we get that
\[\left|\E [X_u^v]-(2n-(x-y)^2n^{1/3})\right|\leq C_2n^{-1/3}(x-y)^4+C_2n^{1/3}\,,\]
for some constants $C_1,C_2>0$, where $X_u^v$ is defined in \eqref{e:defE}. Since $|x-y|\leq 2C_0^{-1}s^{1/4}n^{1/6}$, 
\[C_2n^{-2/3}(x-y)^4\leq 2^4C_0^{-4}C_2s<2^{-1}s\] 
for $C_0>2^{5/4}C_2^{1/4}$. Hence, using the definition of the weight function in \eqref{e:defweight}, for all $s\geq 6C_2$, 
\begin{eqnarray*}
\Bigg\{\left|\weight_{n;(x,0)}^{(y,1)}+(x-y)^2\right|>s\Bigg\}&\subseteq & \Bigg\{n^{-1/3}\left|X_u^v-\E X_u^v\right|>s-C_2n^{-2/3}(x-y)^4-C_2\Bigg\}\\
&\subseteq & \Bigg\{n^{-1/3}\left|X_u^v-\E X_u^v\right|>3^{-1}s\Bigg\}\,.
\end{eqnarray*}

Let $U=T_n^{-1}(I\times \{0\})$ and $V=T_n^{-1}(J\times \{1\})$. For $u\in U,v\in V$, since $3^{-1}<S_{u,v}<3$, we can invoke the proofs of \cite[Propositions~$10.1$ and~$10.5$]{BSS14}. Observe that, for Poissonian last passage percolation, \cite[Corollary~$9.1$]{BSS14} strengthens to 
\begin{equation}\label{e:locc1}
\P(|X_u^{u'}-\E X_u^{u'}|>\theta r^{1/3})\leq e^{-C_1\theta^{3/2}}\,.
\end{equation}
Following the proofs of Proposition $10.1$ and $10.5$ of \cite{BSS14} verbatim, and using the above bound in~\eqref{e:locc1} in place of  Corollary~$9.1$ of \cite{BSS14}, one thus has for all $n,s$ large enough,
\[\P\Bigg(\sup_{u\in U,v\in V}n^{-1/3}\left|X_u^v-\E X_u^v\right|>2^{-1}s\Bigg)\leq e^{-cs^{3/2}}\,.\]
Thus, for $n$ large enough, and $I$ and $J$ intervals of at most unit length contained in the interval of length $2C_0^{-1}s^{1/4}n^{1/6}$ centred at the origin,
\begin{equation}\label{e:01case}
\P \bigg( \, \sup_{x \in I, y \in J} \Big\vert  \weight_{n;(x,0)}^{(y,1)} +  (x-y)^2 \Big\vert > s \, \bigg) \, \leq \,  C \exp \big\{ -c s^{3/2} \big\} \, .
\end{equation}

We now make a first use of the boldface notation for applying results specified at the beginning of Section~\ref{s:estinteg}.
For general $t_1<t_2$, set $\bm{n}=n\tot,\bm{x}=x\tot^{-2/3},\bm{y}=y\tot^{-2/3}, \bm{I}=\tot^{-2/3}I, \bm{J}=\tot^{-2/3}J$ and $\bm{s}=s$ in \eqref{e:01case}.
Recall that the boldface variables are those of Theorem~\ref{t:unifmoddev} and that these are written in terms of non-boldface parameters specified by the present context. 

From the hypothesis of Theorem \ref{t:unifmoddev}, $\bm{I}$ and $\bm{J}$ are intervals of at most unit length contained in $[-\bm{n}^{1/6},\bm{n}^{1/6}]$. 
Thus, applying \eqref{e:01case} and using the scaling principle \eqref{e:eqdist}, we get Theorem~\ref{t:unifmoddev}.
\qed

 The following lower bound on the tail of the polymer weight distribution follows from \cite[Theorem $1.3$]{LM01} and \eqref{e:eqdist}.

\begin{theorem}\label{t:moddevlow}
There exist constants $c_2,s_0,n_0>0$ such that, for all $t_1<t_2$ with $n\tot>n_0$ and $s>s_0$,
\[ \P\left(\tot^{-1/3}\weight_{n;(0,t_1)}^{(0,t_2)}\geq s\right)\geq e^{-c_2s^{3/2}}.\]
\end{theorem}


Moving to unscaled coordinates, the transversal fluctuations for paths between $(0, 0)$ and $(n, n)$  around the interpolating line joining the two points were shown to be $n^{
2/3+o(1)}$ with high
probability in \cite{J00}. More precise estimates were established in \cite{BSS14}. However, the fluctuation of the geodesic at the point $(r,r)$ for any $r\leq n$ is only of the order $r^{2/3}$. This is the content of the next theorem which in essence is the scaled version of \cite[Theorem $2$]{BSS17++} adapted for Poissonian LPP. Recall that, for $u,v\in \R^2$, $\Phi_{n;u}^v$ is the set of all polymers from $u$ to $v$, and $\l_u^v$ is the straight line joining $u$ and $v$. 

\begin{theorem}\label{l:at0}\label{t:carsestimate2} There exist positive constants $n_0, s_1, c$ such that for all $x,y,t_1,t_2\in \R$ with $\tot=t_2-t_1>0$ and $|x-y|\leq 2^{-1}n^{1/3}\tot$ and for all $n\tot\geq n_0, s\geq s_1$ and $t\in [t_1,t_2]$,
\begin{equation}\label{e:cars}
\P\left(\max \left\{ \left|\rho(t)-\l_{(x,t_1)}^{(y,t_2)}(t)\right| : \rho\in \Phi_{n;(x,t_1)}^{(y,t_2)} \right\} \geq s\Big((t-t_1)\wedge (t_2-t)\Big)^{2/3} \right)\leq 2e^{-cs^3}.
\end{equation}
\end{theorem}
Here $a\wedge b$ denotes $\min\{a,b\}$.

\noindent{\bf Proof of Theorem \ref{t:carsestimate2}.} First we prove the theorem when $t_1=0,t_2=1,x=0$. Observe that in this case it is enough to bound the probabilities of the events
$$
\left\{\left|\rho_{n;(0,0)}^{\leftarrow;(y,1)}(t)-\l_{(0,0)}^{(y,1)}(t)\right|\geq s\Big(t\wedge (1-t)\Big)^{2/3}\right\}  \, \, \, \textrm{and} \, \, \, \left\{\left|\rho_{n;(0,0)}^{\rightarrow;(y,1)}(t)-\l_{(0,0)}^{(y,1)}(t)\right|\geq s\Big(t\wedge (1-t)\Big)^{2/3}\right\} \, ,
$$
 and use a union bound to obtain~\eqref{e:cars}.

We first prove an upper bound for the probability of the first of these two events. Also, first assume that $t\in [0,2^{-1}]$. To prove the bound in this case, we move to unscaled coordinates, and use \cite[Theorem $2$]{BSS17++}. 

To this end, let $\Gamma:=\Gamma_{(0,0)}^{\leftarrow;(n+yn^{2/3},n-yn^{2/3})}$  be the leftmost geodesic, and $\mathcal{S}$ the straight line 
from $(0,0)$ to $(n+yn^{2/3},n-yn^{2/3})$. For $r\in [0,n+yn^{2/3}]$, let $\Gamma(r)$ and $\mathcal{S}(r)$ be such that $(r,\Gamma(r))\in \Gamma$ and $(r,\mathcal{S}(r))\in \mathcal{S}$. Now, for $r=nt$, 
\begin{eqnarray}\label{e:scalingfluct}
&&\left\{\left|\rho_{n;(0,0)}^{\leftarrow;(y,\tot)}(t)-\l_{(0,0)}^{(y,\tot)}(t)\right|\geq st^{2/3}\right\}\\
&=&\left\{\left|n^{2/3}\rho_{n;(0,0)}^{\leftarrow;(y,1)}(rn^{-1})-n^{2/3}\l_{(0,0)}^{(y,1)}(rn^{-1})\right|\geq sr^{2/3}\right\}\nonumber\\
&\subseteq &\left\{\left|\Gamma(r')-\mathcal{S}(r')\right|\geq sr^{2/3}\right\}=:\mathsf{B} \,  , \nonumber
\end{eqnarray}
where $r'$ is such that the anti-diagonal line passing through $(r,r)$ intersects $\mathcal{S}$ at $(r',\mathcal{S}(r'))$. The last inclusion follows from the definition of the scaling map $T_n$ in \eqref{e:deftrans}. Since $|y|\leq 2^{-1}n^{1/3}$, $2^{-1}r\leq r'\leq 2r$. Thus,
\begin{equation}
\mathsf{B}\subseteq \left\{ \left|\Gamma(r')-\mathcal{S}(r')\right|\geq 2^{-1}s(r')^{2/3}\right\}=:\mathsf{C}\,. \end{equation}
Thus it is enough to bound the probability of the event $\mathsf{C}$. This local fluctuation estimate for the leftmost geodesic in \eqref{e:local} was proved for exponential directed last passage percolation in \cite[Theorem $2$ and Corollary $2.4$]{BSS17++}. The proof goes through verbatim for the leftmost (and also the rightmost) geodesic in Poissonian last passage percolation. Moreover, the refined bounds of Theorem \ref{t:unifmoddev} give corresponding improvements  for Poissonian LPP: see \cite[Remark $1.5$]{BSS17++}. This gives that, for some positive constants $n_0,r_0,s_0$, and 
for $n\geq n_0, r'\geq r'_0$ and $s\geq s_0$, 
\begin{equation}\label{e:local}
\P(\mathsf{C})\leq e^{-cs^3}\,.
\end{equation}

However, observe that \eqref{e:local} holds only when $r'\geq r'_0$. Now assume $r'\leq r'_0$, so that $r\leq r_0$, where $r_0=2r_0'$. Let the anti-diagonal passing through $(r,r)$ intersect the geodesic $\Gamma$ at $v$ and the line $\mathcal{S}$ at $w$. Clearly $\|v-(r,r)\|_2 \leq 2^{1/2} r$. Also, since $|y|\leq 2^{-1}n^{1/3}$,
\[\|w-(r,r)\|_2= 2^{1/2} |y|rn^{-1/3}\leq r\,.\] 
Thus, with $r=nt\leq r_0$,
\[\left|\rho_{n;(0,0)}^{\leftarrow;(y,1)}(t)-\l_{(0,0)}^{(y,1)}(t)\right|=2^{-1/2}n^{-2/3}\|v-w\|_2 \leq 2^{-1}(2^{1/2}+1)n^{-2/3}r\leq 2r_0^{1/3}t^{2/3}\,. \]
Define $s_1=\max\{s_0,2r_0^{1/3}\}$. Then for $n\geq n_0,s\geq s_1$ and $t\in [0,2^{-1}]$, 
\[\P\left(\left|\rho_{n;(0,0)}^{\leftarrow;(y,1)}(t)-\l_{(0,0)}^{(y,1)}(t)\right|\geq st^{2/3}\right)\leq e^{-cs^3}\,.\]

For $t\in [2^{-1},1]$, we consider the reversed polymer and translate it by $-y$ so that its starting point is $(0,0)$, that is, $\rho'(v)=\rho_{n;(0,0)}^{\leftarrow;(y,1)}(1-v)-y$ for $v\in [0,1]$. Now we follow the same arguments as above to get the bound for the probability of the event 
\[\left\{\left|\rho_{n;(0,0)}^{\leftarrow;(y,1)}(t)-\l_{(0,0)}^{(y,1)}(t)\right|\geq s\Big(t\wedge (1-t)\Big)^{2/3}\right\}\,.\]
Since the same arguments work for the rightmost polymer $\rho_{n;(0,0)}^{\rightarrow;(y,1)}$, we get for $n\geq n_0,s\geq s_1$ and all $t\in [0,1]$,
\begin{equation}\label{e:fix}
\P\left(\max \left\{ \left|\rho(t)-\l_{(0,0)}^{(y,1)}(t)\right| : \rho\in \Phi_{n;(0,0)}^{(y,1)} \right\} \geq s\Big(t\wedge (1-t)\Big)^{2/3} \right)\leq 2e^{-cs^3}.
\end{equation}

Now for general $t_1<t_2$, set $\bm{n}=n\tot, \bm{y}=(y-x)\tot^{-2/3},\bm{s}=s$ and $\bm{t}=\tot^{-1}(t-t_1)$. Then from the hypothesis of Theorem \ref{t:carsestimate2}, $|\bm{y}|\leq 2^{-1}\bm{n}^{1/3}$ since $|y-x|\leq 2^{-1}n^{1/3}\tot$. Thus applying \eqref{e:fix} and using the scaling principle \eqref{e:screl}, we get the theorem.

\qed


The following theorem bounds the transversal fluctuation of polymers; (recall the definitions in \eqref{e:TFfixedpol} and \eqref{e:TFbetwnpoints}). The theorem essentially follows from \cite[Theorem ~$11.1$]{BSS14}; however, we replace the exponent in the upper bound with its optimal value.

\begin{theorem}\label{l:flucpoly}\label{t:transversal}There exist positive constants $c$, $n_0$ and $k_0$  such that, for   $t\in (0,1]$,  $k\geq k_0$ and $n\geq n_0t^{-1}$,
\[\P\left(\mathrm{TF}_{n;(0,0)}^{(0,t)}\geq kt^{2/3}\right)\leq 2e^{-ck^3}.\]
\end{theorem}
\noindent{\bf Proof.} Because of \eqref{e:TFbetwnpoints}, it is enough to bound the probabilities of the events $\left\{\mathrm{TF}\big(\rho_{n;(0,0)}^{\leftarrow;(0,t)}\big)\geq kt^{2/3}\right\}$ and $\left\{\mathrm{TF}\big( \rho_{n;(0,0)}^{\rightarrow;(0,t)} \big) \geq kt^{2/3}\right\}$ and use a union bound. We bound only the first event, the arguments for the second event being the same. Then, as in the proof of Theorem \ref{l:at0}, going to the unscaled coordinates, and defining $\Gamma=\Gamma_{(0,0)}^{\leftarrow;(nt,nt)}$, it is enough to show that
\begin{equation}\label{e:maxTF1}
\P\left(\sup_{r\in [0,nt]}\left|\Gamma(r)-r\right|\geq k(nt)^{2/3}\right)\leq e^{-ck^3}\,.
\end{equation}
From Theorem \ref{t:carsestimate2}, it is easy to see that there exist constants $c>0$ and $n_0,k_0>0$ such that, for all $k>k_0$ and $nt\geq n_0$,
\[\P\left(\left|\Gamma\left(2^{-1}nt\right)-2^{-1}nt\right|>k(nt)^{2/3}\right)\leq e^{-ck^3}\,.\] 
Using the above bound in place of \cite[Lemma $11.3$]{BSS14}, and following the rest of the proof of \cite[Theorem~$11.1$]{BSS14} verbatim, we get \eqref{e:maxTF1}. 
\qed



\section{Basic tools}\label{s:polorder}
Fundamental facts about ordering  and concatenation of polymers will be used repeatedly in  the proofs of the main theorems.

\subsection{Polymer concatenation and superadditivity of weights}\label{ss:pocat}Let $n\in \N$ and $(x,t_1),(y,t_2)\in \R^2$ with $t_1<t_2$ and $|x-y|<n^{1/3}(t_2-t_1)$. (This condition ensures that $(x,t_1)\overset{n}\preceq (y,t_2)$, see Subsection~ \ref{sss:comp}.) Let $u=T_n^{-1}(x,t_1)$ and $v=T_n^{-1}(y,t_2)$ and let $\zeta$ be an increasing path from $u$ to $v$. Let $\gamma=T_n(\zeta)$. We call $\gamma$ an \textit{$n$-path}. We shall often consider $\gamma$ as a subset of $\R^2$, and call $(x,t_1)$ its \textit{starting point} and $(y,t_2)$ its \textit{ending point}. Moreover, similarly to the definition of the weight of a polymer in \eqref{e:defweight}, we define the weight of an $n$-path as
\begin{equation}\label{e:defweightpath}
n^{-1/3}\left(|\zeta|-2n\tot\right)\,,
\end{equation}
where $|\zeta|$ denotes the energy of $\zeta$, that is, the number of points in $\Pi\setminus\{v\}$ that lie on $\zeta$.

Now, let $(x,t_1),(y,t_2),(z,t_3)\in \R^2$ be such that $t_1<t_2<t_3$, $|x-y|< n^{1/3}(t_2-t_1)$ and $|y-z|<n^{1/3}(t_3-t_2)$, so that there exist polymers from $(x,t_1)$ to $(y,t_2)$; and from $(y,t_2)$ to $(z,t_3)$. Let $\rho_1$ be any polymer from $(x,t_1)$ to $(y,t_2)$, and $\rho_2$ any polymer from $(y,t_2)$ to $(z,t_3)$. The union of these two subsets of $\R^2$ is an $n$-path from $(x,t_1)$ to $(z,t_3)$. We call this $n$-path the concatenation of $\rho_1$ and $\rho_2$ and denote it by $\rho_1\circ \rho_2$. The weight of $\rho_1\circ \rho_2$ is $\weight_{n;(x,t_1)}^{(y,t_2)}+\weight_{n;(y,t_2)}^{(z,t_3)}$. This additivity is the reason that the endpoint $v$ was excluded from the definition of path energy in Section~\ref{ss:defres}.

Again, let $n\in \N$ and $(x,t_1),(y,t_2),(z,t_3)\in \R^2$ be such that $t_1<t_2<t_3$ and $|x-y|< n^{1/3}(t_2-t_1)$ and $|y-z|<n^{1/3}(t_3-t_2)$. Then
\begin{equation}\label{e:superadd}
\weight_{n;(x,t_1)}^{(z,t_3)}\geq \weight_{n;(x,t_1)}^{(y,t_2)}+\weight_{n;(y,t_2)}^{(z,t_3)}\,.
\end{equation}
Indeed, taking a polymer $\rho_1$ from $(x,t_1)$ to $(y,t_2)$ and a polymer $\rho_2$ from $(y,t_2)$ to $(z,t_3)$, the weight of $\rho_1\circ\rho_2$ is a lower bound on $\weight_{n;(x,t_1)}^{(z,t_3)}$.

\subsection{Polymer ordering lemmas}
 The first lemma roughly says that if two polymers intersect at two points during their lifetimes, then they are identical between these points. 

\begin{lemma}\label{l:basic}Let $n\in \N$ and $(x_1,t_1),(x_2,t_2),(y_1,s_1),(y_2,s_2)\in \R^2$ and $t,s\in \R$ be such that $t_1<t<s<s_1$, $t_2<t<s<s_2$, $|x_1-y_1|<n^{1/3}(s_1-t_1)$ and $|x_2-y_2|<n^{1/3}(s_2-t_2)$. Suppose that $\rho_{n;(x_1,t_1)}^{\leftarrow;(y_1,s_1)}$ and $\rho_{n;(x_2,t_2)}^{\leftarrow,(y_2,s_2)}$ intersect at two points $z_1=(x,t)$ and $z_2=(y,s)$. Then $\rho_{n;(x_1,t_1)}^{\leftarrow;(y_1,s_1)}$ and $\rho_{n;(x_2,t_2)}^{\leftarrow,(y_2,s_2)}$ are identical between  $t$ and $s$. The same statement holds for the rightmost polymers.
\end{lemma}

To simplify notation in the proof, we write $\rho_1=\rho_{n;(x_1,t_1)}^{\leftarrow;(y_1,s_1)}$ and $\rho_2=\rho_{n;(x_2,t_2)}^{\leftarrow,(y_2,s_2)}$.
\begin{figure}[h]
\includegraphics[scale=.70]{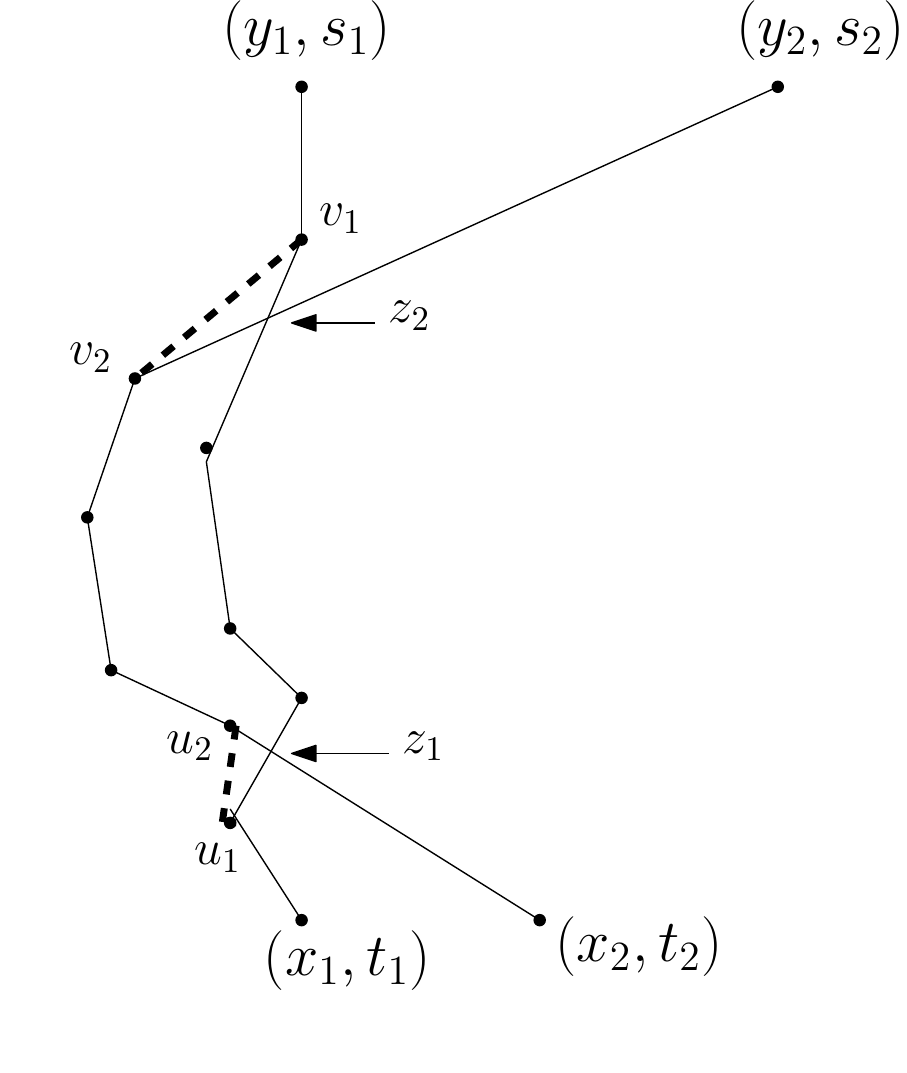}
\caption{This illustrates Lemma \ref{l:porder}. The points of the underlying Poisson process lying on a polymer are marked by dots, and the polymer is obtained by linearly interpolating between the points. The figure shows that both the paths cannot be leftmost polymers between their respective endpoints, since by joining the dashed lines, one obtains an alternative increasing path where the Poisson points between the intersecting points $z_1$ and $z_2$ in the two polymers are interchanged. }
\label{f:PO}
\end{figure}

\noindent{\bf Proof of Lemma~\ref{l:basic}.}
First, for any polymer $\rho$, call a point $u\in \rho$ a \textit{Poisson} point of $\rho$  if $T_n^{-1}(u)\in \Pi\cap\Gamma$, where $\Gamma$ is the geodesic $T_n^{-1}(\rho)$ and $\Pi$ is the underlying unit rate Poisson point process. Also, for $r_1,r_2\in \rho$, let $\rho[r_1,r_2]$ denote the part of the polymer  between the points $r_1$ and~$r_2$, and let $\#\rho[r_1,r_2]$ denote the number of Poisson points that lie in $\rho[r_1,r_2]$. We first claim that $\#\rho_1[z_1,z_2]=\#\rho_2[z_1,z_2|$ where $z_1$ and $z_2$ appear in the lemma's statement. For, if not, without loss of generality assume that $\#\rho_1[z_1,z_2]<\#\rho_2[z_1,z_2|$  and let $u_1$ and $v_1$ be the Poisson points of~$\rho_1$ immediately before $z_1$ and immediately after $z_2$; and let $u_2$ and $v_2$ be the Poisson points of $\rho_2$ immediately after $z_1$ and immediately before $z_2$: see Figure \ref{f:PO}. Then joining $u_1$ to $u_2$ and $v_1$ to $v_2$ (shown in the figure by dashed lines), one gets an alternative path $\rho'$ between $(x_1,t_1)$ and $(y_1,s_1)$ that has more Poisson points than $\rho_1$, thereby contradicting that $\rho_1$ is a polymer between $(x_1,t_1)$ and $(y_1,s_1)$. Thus, $\#\rho_1[z_1,z_2]=\#\rho_2[z_1,z_2|$. Since both $\rho_1$ and $\rho_2$ are leftmost polymers between their respective endpoints, we see that $\rho_1[z_1,z_2]=\rho_2[z_1,z_2]$. This proves the lemma. \qed

The next result roughly says that two polymers that begin and end at the same heights, with the endpoints of one to the right of the other's, cannot cross during their shared lifetime.
\begin{lemma}[Polymer Ordering]
\label{l:porder}
Fix $n\in \N$. Consider points $(x_1,t_1),(x_2,t_1),(y_1,t_2),(y_2,t_2)\in \R^2$ such that $t_1 < t_2$, $x_1\leq x_2$, $y_1\leq y_2$, $|x_1-y_1|<n^{1/3}(t_2-t_1)$ and  $|x_2-y_2|<n^{1/3}(t_2-t_1)$. 
Then $\rho_{n;(x_1,t_1)}^{\leftarrow;(y_1,t_2)}(t)\leq \rho_{n;(x_2,t_1)}^{\leftarrow;(y_2,t_2)}(t)$ and $\rho_{n;(x_1,t_1)}^{\rightarrow;(y_1,t_2)}(t)\leq \rho_{n;(x_2,t_1)}^{\rightarrow;(y_2,t_2)}(t)$ for all $t\in [t_1,t_2]$. 
\end{lemma}
 Let $\rho_1=\rho_{n;(x_1,t_1)}^{\leftarrow;(y_1,t_2)}$ and $\rho_2=\rho_{n;(x_1,t_1)}^{\rightarrow;(y_1,t_2)}$.

\noindent{\bf Proof of Lemma \ref{l:porder}.} Supposing otherwise,  there exists $z=(x,y)\in \rho_2$ such that $x<\rho_1(y)$. But then there exist $z_1,z_2\in \rho_1 \cap\rho_2$ straddling the point $z$. By Lemma \ref{l:basic}, $\rho_1[z_1,z_2]=\rho_2[z_1,z_2]$, and hence $z\in \rho_1\cap \rho_2$, a contradiction.
\qed

By ordering, a polymer whose endpoints are straddled between those of a pair of polymers becomes sandwiched between those polymers.

\begin{corollary}\label{c:po}Fix $n\in \N$. Consider points $(x_1,t_1),(x_2,t_1),(x_3,t_1),(y_1,t_2),(y_2,t_2),(y_3,t_2)\in \R^2$ such that $t_1 < t_2$, $x_1\leq x_2\leq x_3$, $y_1\leq y_2\leq y_3$ and $|x_i-y_i|<n^{1/3}(t_2-t_1)$ for $i=1,2,3$. Let $t\in (t_1,t_2)$. Let $\rho_i=\rho_{n;(x_i,t_1)}^{\leftarrow;(y_i,t_2)}$ for $i=1,2,3$. Then 
\[|\rho_2(t)-\rho_2(t_1)|\leq \max_{i\in\{1,3\}}|\rho_i(t)-\rho_i(t_1)|+\max_{i\in\{1,3\}}|x_i-x_2|\,.\]
The same result holds for rightmost polymers.
\end{corollary}
\noindent{\bf Proof.} By Lemma \ref{l:porder},
\[\rho_1(t)\leq \rho_2(t)\leq \rho_3(t)\,.\]
The result now follows immediately.
\qed



\section{Exponent pair $(2/3,1/3)$ for a single polymer: Proof of Theorem \ref{t:main1}}\label{s:2/3,1/3}
In this section, we show that the sequence 
$\big\{\rho^{\leftarrow}_n : n \in \N \big\}$ of leftmost $n$-polymers from $(0,0)$ to $(0,1)$ is tight, and any weak limit is \holder $2/3-$-continuous with a polylogarithmic correction of order~$1/3$. The main two ingredients in this proof are the local regularity estimate Theorem \ref{t:carsestimate2} and the polymer ordering Lemma \ref{l:porder}. First, we bound the fluctuation of the polymer near any given point $z\in [0,1]$.

\begin{proposition}\label{l:boundatz}
There exist positive constants $n_0,s_1$ and $c$  such that, for all $n\geq n_0,s\geq s_1$, $z\in [0,1]$ and $0\leq t\leq 1-z$,
\begin{equation}\label{e:boundatz}
\P \left( |\rho_n^{\leftarrow}(z+t)-\rho_n^{\leftarrow}(z)|\geq st^{2/3} \right) \leq 10t^{-2/3}e^{-cs^3}.
\end{equation}
The same statement holds for $\rho_n^{\rightarrow}$.
\end{proposition}

\begin{figure}[h]
\centering
\includegraphics[width=0.8\textwidth]{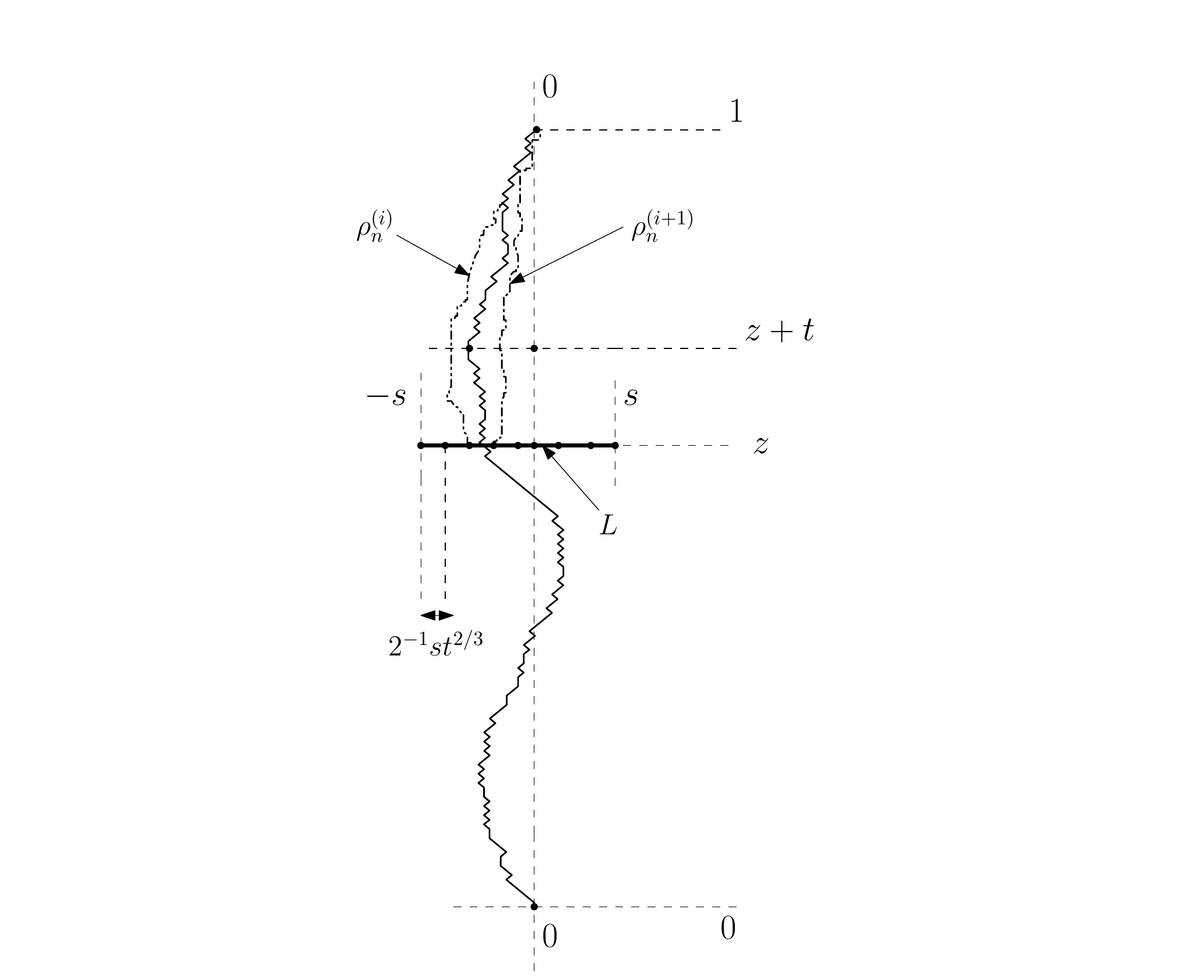}
\caption{The proof of Proposition \ref{l:boundatz} is illustrated here. We mark the line segment $L$ with a number of equally spaced points. As the leftmost polymer from $(0,0)$ to $(0,1)$ passes between two such points on the line $L$, it is, in view of polymer ordering, sandwiched between the two leftmost polymers, shown as dotted lines, originating from those points and ending at $(0,1)$. Hence it is sufficient to bound the fluctuations of the polymers originating from these equally spaced points on $L$.}
\label{f:Fig1}
\end{figure}

As we now explain, the proposition will be proved by reducing to the case that $z=0$, when the result  follows from Theorem \ref{l:at0}. For any fixed $z\in (0,1)$, Theorem \ref{l:at0} again guarantees that the polymer $\rho_n^{\leftarrow}$ is at distance at most $s$ from the point $(0,z)$ with probability at least $1-e^{-cs^3}$. We break the horizontal line segment of length $2s$ centred at $(0,z)$ into a sequence of consecutive intervals of length $2^{-1}st^{2/3}$, and consider the leftmost polymers starting from each of these endpoints and ending at $(0,1)$, as in Figure \ref{f:Fig1}. Due to the Corollary \ref{c:po} of  the polymer ordering Lemma \ref{l:porder}, a big fluctuation of $\rho_n^{\leftarrow}$ between times $z$ and $z+t$ creates a big fluctuation for one of the polymers starting from these deterministic endpoints. The probability of the latter  fluctuations is controlled via Theorem \ref{l:at0} and since the number of these polymers is of the order of $t^{-2/3}$, a union bound gives~\eqref{e:boundatz}. \ \\
 
\noindent{\bf Proof of Proposition \ref{l:boundatz}.} 
First observe that for $s> (nt)^{1/3}$,  the  probability in \eqref{e:boundatz} is zero by the definition of the scaling map $T_n$ in \eqref{e:deftrans} and the geodesics being increasing paths. Hence we assume that $s\leq (nt)^{1/3}$. 

Fix $s\leq (nt)^{1/3}$ and $z\in [0,1]$. For $t\geq 8^{-3}$, 
\[\left\{|\rho_n^{\leftarrow}(z+t)-\rho_n^{\leftarrow}(z)|\geq st^{2/3}\right\}\subseteq \left\{|\rho_n^{\leftarrow}(z+t)-\rho_n^{\leftarrow}(z)|\geq 8^{-2}s\right\}\subseteq\left\{\mathrm{TF}_{n;(0,0)}^{(0,1)}\geq 2^{-1}8^{-2}s\right\}\,,\]
where $\mathrm{TF}_{n;(0,0)}^{(0,1)}$ is defined in \eqref{e:TFbetwnpoints}. Hence, applying Theorem \ref{l:flucpoly} with the parameter specifications $\bm{t}=1$ and $\bm{k}=2^{-1}8^{-2}s$, we get that \eqref{e:boundatz} holds for all $n,s$ large enough. Hence we assume that $t\leq 8^{-3}$. Also, let us assume for now that $z\in [0,2^{-1}]$. 

Let $L$ be the line segment $[-s,s]\times\{z\}$. Let $\mathsf{E}$ be the event that $\rho_n^{\leftarrow}$ passes through~$L$. By Theorem~\ref{l:at0} with $\bm{n}=n, \bm{t}=z, \bm{x}=0,\bm{y}=0, \bm{s}=s,\bm{t_1}=0$ and $\bm{t_2}=1$, we have that, for $n\geq n_1$ and $s\geq s_1$, 
\[\P(\mathsf E)\geq 1-2e^{-cs^3}.\]


Now, we divide $L$ into $\lceil{4t^{-2/3}}\rceil$-many adjacent intervals of length at most ${2^{-1}st^{2/3}}$, and let $(x_i,z), i=0,1,2,\cdots,\lceil4t^{-2/3}\rceil$ be the endpoints of these intervals, i.e., 
\[x_i=-s+2^{-1}ist^{2/3} \quad \mbox{for } i=0,1,2,\cdots, \lceil4t^{-2/3}\rceil.\]

Let $\rho_n^{(i)}:=\rho_{n;(x_i,z)}^{\leftarrow;(0,1)}$ be the leftmost polymer from $(x_i,z)$ to $(0,1)$. 

By Corollary \ref{c:po}, on $\mathsf{E}$,
\begin{align}\label{e:polord}
\left|\rho_n^{\leftarrow}(z+t)-\rho_n^{\leftarrow}(z)\right|\leq \max_{i\in \llbracket 0, \lceil4t^{-2/3}\rceil\rrbracket}\left|\rho_n^{(i)}(z+t)-\rho_n^{(i)}(z)\right|
+2^{-1}st^{2/3}.
\end{align}

Also, for any fixed $i \in \llbracket 0, \lceil4t^{-2/3}\rceil\rrbracket$, let $\l^{(i)}=\l_{(x_i,z)}^{(0,1)}$ be the straight line segment joining $(x_i,z)$ and $(0,1)$. Then, since $z\in [0,2^{-1}]$ and $t\leq 8^{-3}$, for any $i\in 0,1,2,\cdots,\lceil4t^{-2/3}\rceil$,
$$\left|\l^{(i)}(z)-\l^{(i)}(z+t)\right|
\leq
\frac{st}{1-z}
\leq 
2st
\leq 
4^{-1}st^{2/3} \, .
$$

Since $\rho_n^{(i)}(z)=\l^{(i)}(z)=x_i$,
\begin{eqnarray*}
\left|\rho_n^{(i)}(z+t)-\rho_n^{(i)}(z)\right|
&\leq& \left|\rho_n^{(i)}(z+t)-\l^{(i)}(z+t)\right|+|\l^{(i)}(z+t)-\l^{(i)}(z)|\\
&\leq& \left|\rho_n^{(i)}(z+t)-\l^{(i)}(z+t)\right|+{4^{-1}st^{2/3}}\,.
\end{eqnarray*}
Thus, on the event $\mathsf{E}$, by \eqref{e:polord},
\begin{equation*}
\left|\rho_n^{\leftarrow}(z+t)-\rho_n^{\leftarrow}(z)\right|\leq \max_{i\in \llbracket 0, \lceil4t^{-2/3}\rceil\rrbracket}\left|\rho_n^{(i)}(z+t)-\l^{(i)}(z+t)\right|+ 
\tfrac{3}{4} st^{2/3} \, .
\end{equation*}

From here, it follows by taking a union bound that
\begin{eqnarray*}
&&\P\left(\left|\rho_n^{\leftarrow}(z+t)-\rho_n^{\leftarrow}(z)\right|\geq st^{2/3}\right)\\
&\leq& \P(\mathsf{E}^c)+\sum_{i=0}^{\lceil 4t^{-2/3}\rceil}\P\left(\left|\rho_n^{(i)}(z+t)-\l^{(i)}(z+t)\right|\geq  4^{-1}st^{2/3} \right)\\
&\leq & 10t^{-2/3}e^{-cs^3},
\end{eqnarray*}
for some absolute positive constant $c$ and all $n\geq 2n_0$. Here the last inequality follows by applying Theorem~\ref{l:at0} to each of the polymers $\rho^{(i)}$. For given $i$, set the parameters $\bm{n}=n,\bm{t_1}=z,\bm{t_2}=1,\bm{t}=t+z,\bm{x}= -s+2^{-1}ist^{2/3},\bm{y}=0$ and $\bm{s}=4^{-1}s$. Since $z\in [0, 2^{-1}]$ and $s\leq (nt)^{1/3}$, we have that  $|\bm{x}-\bm{y}|\leq s\leq n^{1/3}t^{1/3}\leq 8^{-1}n^{1/3}\leq 4^{-1}\bf{n}^{1/3}\bf{\tot}$. Thus one can apply Theorem~\ref{l:at0} to get the above inequality for all $\bm{n\tot}\geq 2^{-1}n\geq n_0$.

For $z\in [2^{-1},1]$, define the reversed polymer $\widehat{\rho}_n^{\leftarrow}$ by $\widehat{\rho}_n^{\leftarrow}(a)={\rho}_n^{\leftarrow}(1-a)$ for $a\in [0,1]$, and follow the above argument. 
\qed

Next we show the tightness of the members of the sequence $\{\rho_n^{\leftarrow}\}_{n\in \N}$ as elements in the space $(C[0,1],\|\cdot\|_\infty)$. We prove that Proposition \ref{l:boundatz} guarantees that Kolmogorov-Chentsov's tightness criterion is satisfied.

\noindent{\bf Proof of Theorem \ref{t:main1}$(a)$.} Fix $n\geq n_0$ and any $\lambda>0$. Fix $t\in (0,1]$ small enough  that $\lambda t^{-2/3}\geq s_1$, where $n_0$ and $s_1$ are as in Proposition \ref{l:boundatz}. Also fix some $M\in \N$ large enough  that $2M-2/3>1$. Then it follows from Proposition \ref{l:boundatz} that for any $z,z'\in [0,1]$ with $|z-z'|=t$,
\begin{equation}\label{e:KolCh}
\P\left(|\rho_n^{\leftarrow}(z')-\rho_n^{\leftarrow}(z)|\geq \lambda\right)\leq 10t^{-2/3}e^{-c(\lambda^3 t^{-2})}\leq K_M\lambda^{-3M} t^{2M-2/3}=K_M\lambda^{-3M} |z'-z|^{2M-2/3},
\end{equation}
where $K_M:=\sup_{x\geq 0}x^Me^{-cx}<\infty$. Since $2M-2/3>1$, by Kolmogorov-Chentsov's tightness criterion (see for example \cite[Theorem~$8.1.3$]{Dur10}), it follows that the sequence $\{\rho_n^{\leftarrow}\}_{n\in \N}$ is tight in $(C[0,1],\|\cdot\|_\infty)$.
\qed

\subsection{Modulus of continuity}
 Here we prove Theorem \ref{t:main1}(b), thus finding  the modulus of continuity for any weak limit of a weakly converging subsequence of $\{\rho_n^{\leftarrow}\}_{n\in \N}$. We will follow the arguments used to derive the Kolmogorov continuity criterion, where one infers \holder continuity of a stochastic process from moment bounds on the difference of the process between pairs of times.  Thus we introduce  the set of dyadic rationals
\[D=\bigcup_{i=0}^\infty 2^{-i}\Z \, .\]
Next is the first step towards proving the modulus of continuity. 
\begin{lemma}\label{l:dyad} Let $\rho_*^{\leftarrow}$ be the weak limit of a weakly converging subsequence of $\{\rho_n^{\leftarrow}\}_{n\in \N}$. Then there exists a universal positive constant $C$ (not depending on the particular weak limit $\rho_*^{\leftarrow}$) such that, almost surely, for some random $m_0(\omega)\in \N$ and for all $s,t\in D\cap [0,1]$ with $|t-s|\leq 2^{-m_0(\omega)}$,
\[|\rho_*^{\leftarrow}(t)-\rho_*^{\leftarrow}(s)|\leq C(t-s)^{2/3}\left(\log({t-s})^{-1}\right)^{1/3} .\]
\end{lemma}
\noindent{\bf Proof.} 
For $m\in \N$, let $S_m$ be the set of all intervals of the form $[j2^{-m},(j+1)2^{-m}]$,
for $j\in\{0,1,2,\cdots,2^m-1\}$.
Fix $c_0>(\frac{5}{3c})^{1/3}$, where $c$ is the constant in Proposition \ref{l:boundatz}. 

Writing  $\Rightarrow$ for convergence in distribution, let $\{\rho_{n_k}^{\leftarrow}\}_{k\in \N}$ be a subsequence of $\{\rho^{\leftarrow}_n\}_{n\in \N}$ such that $\rho_{n_k}^{\leftarrow}\Rightarrow \rho_*^{\leftarrow}$ as random variables in $(C[0,1],\|\cdot\|_\infty)$. 
 Since for $a,b\in [0,1]$, the map $\tau_{a,b}$ defined by $(C[0,1],\|\cdot\|_\infty)\mapsto (\R,|\cdot|):f\mapsto |f(a)-f(b)|$ is continuous, 
  \[\mathsf{U}:=\bigcup\left\{\tau^{-1}_{(j+1)2^{-m},j2^{-m}}\left(c_02^{-\frac{2m}{3}}\left(\log 2^m\right)^{1/3},\infty\right):j=0,1,\cdots,2^m-1\right\}\]
 is an open set. Thus, by the Portmanteau theorem,
 \begin{eqnarray}\label{e:main1}
&&\P\Bigg(\sup_{j\in \{0,1,\cdots,2^m-1\}}|\rho_*^{\leftarrow}((j+1)2^{-m})-\rho_*^{\leftarrow}(j2^{-m})|>c_02^{-\frac{2m}{3}}\left(\log 2^m\right)^{1/3}\Bigg)\\
&\leq & \liminf_k \, \P\Bigg(\sup_{j\in \{0,1,\cdots,2^m-1\}}|\rho_{n_k}^{\leftarrow}((j+1)2^{-m})-\rho_{n_k}^{\leftarrow}(j2^{-m})|>c_02^{-\frac{2m}{3}}\left(\log 2^m\right)^{1/3}\Bigg)\nonumber\\
&\leq & \limsup _n \, \P\Bigg(\sup_{j\in \{0,1,\cdots,2^m-1\}}|\rho_{n}^{\leftarrow}((j+1)2^{-m})-\rho_{n}^{\leftarrow}(j2^{-m})|>c_02^{-\frac{2m}{3}}\left(\log 2^m\right)^{1/3}\Bigg)\nonumber\,.
\end{eqnarray}
Now, for all $m$ large enough that $\left(\log 2^m\right)^{1/3}\geq s_1$, where $s_1$ is as in Proposition \ref{l:boundatz}, and all $n\geq n_0$, applying Proposition \ref{l:boundatz} and a union bound, 
\begin{eqnarray*}
&&\P \, \Bigg(\sup_{j\in \{0,1,\cdots,2^m-1\}}|\rho_n^{\leftarrow}((j+1)2^{-m})-\rho_n^{\leftarrow}(j2^{-m})|>c_02^{-\frac{2m}{3}}\left(\log 2^m\right)^{1/3}\Bigg)\nonumber\\
&\leq & 10\cdot 2^m\left(\frac{1}{2^{m}}\right)^{c_0^3c-2/3}\leq 10\left(\frac{1}{2^{m}}\right)^{c_0^3c-5/3}.\nonumber
\end{eqnarray*}
Hence, from \eqref{e:main1},
\[\P \, \Bigg(\sup_{j\in \{0,1,\cdots,2^m-1\}}|\rho_*^{\leftarrow}((j+1)2^{-m})-\rho_*^{\leftarrow}(j2^{-m})|>c_02^{-\frac{2m}{3}}\left(\log 2^m\right)^{1/3}\Bigg)\leq 10\cdot 2^{-m(c_0^3c-5/3)}\,.\]
As the right hand side is summable in $m$ (by the choice of $c_0$ made at the beginning of the proof), the Borel-Cantelli lemma implies that there exists a null set $N_0$, such that, for each $\omega \notin N_0$, there is some $m_0(\omega)$ for which $m\geq m_0(\omega)$ entails that
\begin{equation}\label{e:dyad}
|\rho_*^{\leftarrow}(t)-\rho_*^{\leftarrow}(s)|\leq c_0(t-s)^{2/3}\left(\log({t-s})^{-1}\right)^{1/3} \quad \mbox{ for all } [s,t]\in S_m \, .
\end{equation}
Now, let $\omega \notin N_0$ and $s,t\in D\cap [0,1]$ be such that $|s-t|\leq 2^{-m_0(\omega)}$.
Let $m=m(s,t)$ be the greatest integer such that $|s-t|\leq 2^{-m}$; then clearly, $m\geq m_0(\omega)$.
Also, consider the binary expansions of $s$ and $t$:
\[s=s_0+\sum_{j>m} \sigma_j2^{-j}, \quad t=t_0+\sum_{j>m} \tau_j2^{-j},\]
where $\sigma_j,\tau_j\in \{0,1\}$, and each of the sequences is eventually zero. Either $s_0=t_0$ or $[s_0,t_0]\in S_m$. Moreover, for $n\geq 1$, let
\[s_n=s_0+\sum_{m<j\leq m+n} \sigma_j2^{-j}.\]
Then, for $n\geq 1$, either $s_n=s_{n-1}$ or $[s_{n-1},s_n]\in S_{m+n}$. Since $m\geq m_0(\omega)$, by \eqref{e:dyad},
\[|\rho_*^{\leftarrow}(t_0)(\omega)-\rho_*^{\leftarrow}(s_0)(\omega)|\leq c_02^{-\frac{2m}{3}}\left(\log 2^{m}\right)^{1/3}.\]
Also,
\begin{eqnarray*}
|\rho_*^{\leftarrow}(s)(\omega)-\rho_*^{\leftarrow}(s_0)(\omega)|\leq \sum_{n=1}^\infty|\rho_*^{\leftarrow}(s_n)(\omega)-\rho_*^{\leftarrow}(s_{n-1})(\omega)|&\leq & \sum_{n=1}^\infty c_02^{-\frac{2(m+n)}{3}}\left(\log 2^{m+n}\right)^{1/3}\\
&\leq & C_12^{-\frac{2(m+1)}{3}}\left(\log 2^{m+1}\right)^{1/3},
\end{eqnarray*}
and similarly
\[|\rho_*^{\leftarrow}(t)(\omega)-\rho_*^{\leftarrow}(t_0)(\omega)|\leq  C_2 2^{-\frac{2(m+1)}{3}}\left(\log 2^{m+1}\right)^{1/3},\]
for some absolute constants $C_1$ and $C_2$. Hence,
\[|\rho_*^{\leftarrow}(t)-\rho_*^{\leftarrow}(s)|\leq |\rho_*^{\leftarrow}(t)-\rho_*^{\leftarrow}(t_0)|+|\rho_*^{\leftarrow}(t_0)-\rho_*^{\leftarrow}(s_0)|+|\rho_*^{\leftarrow}(s)-\rho_*^{\leftarrow}(s_0)|\leq C2^{-\frac{2m}{3}}\left(\log 2^{m}\right)^{1/3}.\]
Since by definition $2^{-m-1}\leq |s-t|\leq 2^{-m}$, the result follows. 
\qed

\noindent{\bf Proof of Theorem \ref{t:main1}$(b)$.}
 For any $s,t\in [0,1]$ satisfying $s<t$ and $|s-t|\leq 2^{-m_0(\omega)}$, choose $s_k,t_k\in D\cap [s,t]$ such that $s_k\searrow s$ and $t_k\nearrow t$. Then, since $|s_k-t_k|\leq |s-t|\leq 2^{-m_0(\omega)}$, by Lemma~\ref{l:dyad}, 
\[|\rho_*^{\leftarrow}(t_k)-\rho_*^{\leftarrow}(s_k)|\leq C(t_k-s_k)^{2/3}\left(\log({t_k-s_k})^{-1}\right)^{1/3} .\]
Since $\rho_*^{\leftarrow}(t_k)(\omega)\rightarrow \rho_*^{\leftarrow}(\omega)$ and $\rho_*^{\leftarrow}(s_k)(\omega)\rightarrow \rho_*^{\leftarrow}(s)(\omega)$, the theorem follows by taking the limit as $k\rightarrow \infty$. The same argument applies without any change for the rightmost polymers as well. 
\qed

\section{Exponent pair $(2/3,1/3)$ for maximum fluctuation over short polymers: \\
Proof of Theorem \ref{t:flucshortpol} }\label{s:shortpol}

In this section, we shall prove Theorem \ref{t:flucshortpol}. It is the upper bound that is the more subtle. Recall the notation of transversal fluctuations from \eqref{e:TFfixedpol} and \eqref{e:TFbetwnpoints}, $\mathrm{AdEndPair}_n(t)$ from \eqref{e:defComp} and $\mathrm{MTF}_n(t)$ from~\eqref{e:defMTF}. 

Here is the idea behind the proof.
Proposition~\ref{p:lowbndTF} offers a lower bound on the transversal fluctuation of a polymer between two given points. By considering order-$t^{-1}$ endpoint pairs with disjoint intervening lifetimes of length $t$,
we obtain a collection of independent opportunities for the fluctuation lower bound to occur. By tuning the probability of the individual event to have order $t$, at least one among the constituent events typically does occur, and the lower bound in Theorem~\ref{t:flucshortpol} follows. 

\begin{figure}
\centering

\includegraphics[scale=.6]{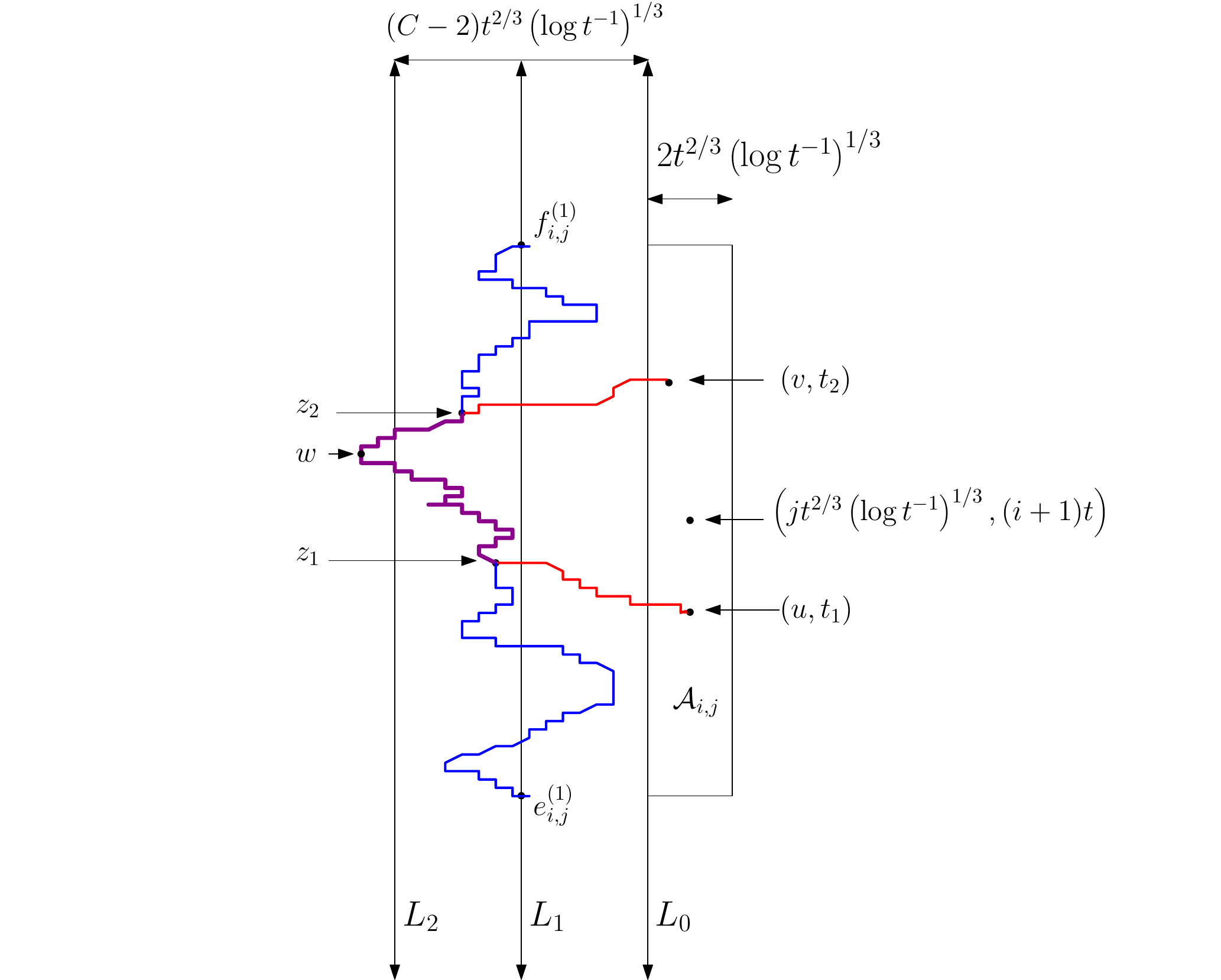} 
\caption{The figure illustrates the proof of the upper bound in Theorem \ref{t:flucshortpol}. If the leftmost polymer between $(u,t_1)$ and $(v,t_2)$ (shown in red) makes a huge leftward fluctuation and the leftmost polymer between points $e_{i,j}^{(1)}$ and $f_{i,j}^{(1)}$ (shown in blue) is to the left of $u$ and $v$ at $t_1$ and $t_2$ respectively, then the blue polymer stays to the left of the red polymer between times $t_1$ and $t_2$ by polymer ordering. Thus the big left fluctuation transmits from the red to the blue polymer. If, however, the blue polymer reaches to the right of either $u$ or $v$, then it creates a big right fluctuation for the blue polymer. Thus by bounding the fluctuations of a small number of polymers between deterministic endpoints, one can bound the fluctuation between all admissible endpoint pairs.}
\label{f:shortpoly}
\end{figure}

On the other hand, suppose that a big swing in the unit order region happens between a certain endpoint pair,
with an intervening duration, or height difference, of order $t$. Members of the endpoint pair may be exceptional locations when viewed as functions of the underlying Poisson point field, both in horizontal and vertical coordinate. 
Thus, the upper bound in Theorem~\ref{t:flucshortpol}
does not follow directly from a union bound of a given endpoint estimate over elements in a discrete mesh, since such a mesh may not capture the exceptional endpoints. However, polymer ordering forces exceptional behaviour to become typical and to occur between an endpoint pair in a discrete mesh. 
To see this, assume that the original polymer between exceptional endpoints makes a big left swing. (Figure \ref{f:shortpoly} illustrates the argument.) 
We take a discrete mesh endpoint pair 
whose lifetime includes that of the original polymer but has the same order $t$, and whose lower and upper points lie to the left of the original endpoint locations, about halfway between these and the leftmost coordinate visited by the original polymer. Then we consider the leftmost mesh polymer at the beginning and ending times of the original polymer. If the mesh polymer is to the right of the original polymer at any of these endpoints, then the mesh polymer has already made a big rightward swing at one of these endpoints. If, on the other hand, the mesh polymer is to the left of the original polymer at both the endpoints of the original polymer, then by polymer ordering Lemma \ref{l:porder}, the mesh polymer cannot cross the original polymer during the latter's lifetime. Hence the big left swing of the original polymer forces a significant left swing for the mesh polymer as well.



\ \\
\noindent{\bf Proof of Theorem \ref{t:flucshortpol}.}
The lower bound follows in a straightforward way from Proposition~\ref{p:lowbndTF}. For any $t\in (0,1)$ and $i\in \left\{0,1,2,\cdots,\left[t^{-1}\right]-1\right\}$, define
\[\mathsf{F}_{i,t,n}=\left\{\mathrm{TF}_{n;(0,it)}^{(0,(i+1)t}\geq ct^{2/3}\big(\log t^{-1}\big)^{1/3}\right\}.\]
For given such $(t,i)$, we  apply Proposition \ref{p:lowbndTF} with parameter settings $\bm{n}=n,\bm{t_1}=it,\bm{t_2}=(i+1)t$ and $\bm{s}=c(\log t^{-1}\big)^{1/3}$,  to find that, when $c(\log t^{-1})^{1/3}\geq s_0$ and $n\geq \max\{\alpha_0^{-3}c^3t^{-1}\log t^{-1},n_0t^{-1}\}$,
 \[\P(\mathsf{F}_{i,t,n})\geq e^{-c^*c^3\log t^{-1}}=t^{c^*c^3 }\,,\]
 where the proposition specifies the quantities $\alpha_0,n_0$ and $s_0$.
 
 Thus, for all $t\leq e^{-(c^{-1}s_0)^3}$ and $i\in \left\{0,1,2,\cdots,\left[t^{-1}\right]-1\right\}$,
\begin{eqnarray*}
\liminf_n t^{-1}\P(\mathsf{F}_{i,t,n})=\liminf_n t^{-1}\P(\mathsf{F}_{0,t,n})
\geq t^{c^*c^3 - 1}\,.
\end{eqnarray*}

By choosing $c>0$ small enough that $c^*c^3<1$, one has
$\liminf_n t^{-1}\P(\mathsf{F}_{0,t,n})\rightarrow \infty$ as $t \searrow 0$.
For such $c>0$, using the definition~(\ref{e:defMTF}) of $\mathrm{MTF}_n(t)$ and independence of the events $\mathsf{F}_{i,t,n}$ for $i\in \left\{0,1,2,\cdots,\left[t^{-1}\right]-1\right\}$,
\[\P\left({\mathrm{MTF}_n(t)}{t^{-2/3}\big(\log t^{-1}\big)^{-1/3}}<c\right) \, \leq \,  \P\left(\bigcap_{i=0}^{[t^{-1}]-1}\mathsf{F}_{i,t,n}^c\right)=\prod_{i=0}^{[t^{-1}]-1}\P\left(\mathsf{F}_{i,t,n}^c\right)\,.\]
Thus,
\begin{eqnarray*}
&&\limsup_n \, \P\left({\mathrm{MTF}_n(t)}{t^{-2/3}\big(\log t^{-1}\big)^{-1/3}}<c\right)\\
&\leq&\limsup_n  \, \Big(1-\P\big(\mathsf{F}_{0,t,n}\big) \Big)^{[t^{-1}]}\leq  \limsup_n  \,  \exp\left\{ -\left[t^{-1}\right]\P(\mathsf{F}_{0,n})\right\} \, \rightarrow \,  0\,,
\end{eqnarray*}
the latter convergence as $t\searrow 0$.

Now we show the upper bound. 
Fix $t\in (0,1]$ small enough that $\psi t\leq t^{2/3}$, where the parameter~$\psi$ appears in the definition~\eqref{e:defComp} of  $\mathrm{AdEndPair}_n(t)$.


For any $i=0,1,2,\ldots, \lceil t^{-1}\rceil$ and $j\in \left\llbracket-\left\lceil t^{-2/3}\big(\log t^{-1}\big)^{-1/3}\right\rceil, \left\lceil t^{-2/3}\big(\log t^{-1}\big)^{-1/3}\right\rceil \right\rrbracket$, define the rectangle $\mathcal{A}_{i,j}$ with lower-left corner $\left((j-1)t^{2/3}\big(\log t^{-1}\big)^{1/3},it\right)$,
width $2 t^{2/3}\big(\log t^{-1}\big)^{1/3}$ and height~$2t$.
Figure \ref{f:shortpoly} illustrates this rectangle and the arguments that follow.

Let $C > 0$ be an even integer whose value will later be specified. 
 For such $i,j$ as above, define planar points 
 \begin{eqnarray*}
 e_{i,j}^{(1)}&:=&\left((j-2^{-1}C)t^{2/3}\big(\log t^{-1}\big)^{1/3},it\right),\quad f_{i,j}^{(1)}:=\left((j-2^{-1}C)t^{2/3}\big(\log t^{-1}\big)^{1/3},(i+2)t\right),\\
 e_{i,j}^{(2)}&:=&\left((j+2^{-1}C)t^{2/3}\big(\log t^{-1}\big)^{1/3},it\right),\quad f_{i,j}^{(2)}:=\left((j+2^{-1}C)t^{2/3}\big(\log t^{-1}\big)^{1/3},(i+2)t\right).
\end{eqnarray*}

Then we claim that, whatever the value of $C > 0$,
\begin{equation}\label{e:imp}
\mathsf{B}_{i,j}:=\left\{\sup\left\{\mathrm{TF}_{n;(x_1,y_1)}^{(x_2,y_2)}:(x_1,y_1),(x_2,y_2)\in \mathcal{A}_{i,j},y_2>y_1\right\}>Ct^{2/3}\big(\log t^{-1}\big)^{1/3}\right\}\subseteq \mathsf{D}^{(1)}_{i,j}\cup \mathsf{D}^{(2)}_{i,j},
\end{equation}
where  \[\mathsf{D}^{(1)}_{i,j}:=\left\{\mathrm{TF}_{n;e_{i,j}^{(1)}}^{f_{i,j}^{(1)}}\ge(2^{-1}C-1)t^{2/3}\big(\log t^{-1}\big)^{1/3}\right\}\]
and
\[\mathsf{D}^{(2)}_{i,j}:=\left\{\mathrm{TF}_{n;e_{i,j}^{(2)}}^{f_{i,j}^{(2)}}\ge(2^{-1}C-1)t^{2/3}\big(\log t^{-1}\big)^{1/3}\right\}.\]
To see \eqref{e:imp}, define the vertical lines: 
$$
L_2 =  \big\{ x = (j-C+1)t^{2/3}\big(\log t^{-1}\big)^{1/3} \big\}  \, \, \, \, \textrm{and} \, \, \, \,
 L_2' = \big\{ x = (j+C-1)t^{2/3}\big(\log t^{-1}\big)^{1/3} \big\} \, .
$$
 Then, on the event $\mathsf{B}_{i,j}$, there exists a pair of points $(u,t_1),(v,t_2)\in \mathcal{A}_{i,j}$ such that either  $\rho_{n;(u,t_1)}^{\leftarrow;(v,t_2)}$ intersects $L_2$ or $\rho_{n;(u,t_1)}^{\rightarrow;(v,t_2)}$ intersects $L_2'$. We now show that, when $\rho_{n;(u,t_1)}^{\leftarrow;(v,t_2)}$ intersects $L_2$, 
 the event $\mathsf{D}_{i,j}^{(1)}$ occurs. 
   Let \[\rho:=\rho_{n;e_{i,j}^{(1)}}^{\leftarrow;f_{i,j}^{(1)}}\,.\]
Let $\l_{i,j}^{(1)}$ be the line segment joining $e_{i,j}^{(1)}$ and $f^{(1)}_{i,j}$. If $\rho(t_1)>u$, then 
\[\rho(t_1)-\l_{i,j}^{(1)}(t_1)\geq (j-1)t^{2/3}\big(\log t^{-1}\big)^{1/3}-(j-2^{-1}C)t^{2/3}\big(\log t^{-1}\big)^{1/3}\geq (2^{-1}C-1)t^{2/3}\big(\log t^{-1}\big)^{1/3}\,,\]
and thus $\mathsf{D}_{i,j}^{(1)}$ holds. Similarly, if $\rho(t_2)>v$, then $\mathsf{D}_{i,j}^{(1)}$ holds. Now assume that $\rho(t_1)<u$ and $\rho(t_2)<v$. Polymer ordering Lemma~\ref{l:porder} then implies that $\rho(t)\leq \rho_{n;(u,t_1)}^{\leftarrow;(v,t_2)}(t)$ for all $t\in [t_1,t_2]$. Thus $\rho$ intersects $L_2$ as well, and hence $\mathsf{D}_{i,j}^{(1)}$ occurs. 
 
  By similar reasoning, we see that, when $\rho_{n;(u,t_1)}^{\rightarrow;(v,t_2)}$ intersects $L_2'$, the event  $\mathsf{D}_{i,j}^{(2)}$ occurs. We have proved~\eqref{e:imp}. 
  
 For any compatible pair of points $(u,v)\in \mathrm{AdEndPair}_n(t)$, there exists a pair $(i,j)$ for which $u,v\in \mathcal{A}_{i,j}$; here we use $\psi t\leq t^{2/3}$. 
Hence,
\begin{eqnarray*}
 & & \left\{{  t^{-2/3}\big(\log t^{-1}\big)^{-1/3} \mathrm{MTF}_n(t)}    >  C\right\}\\
&\subseteq & \bigcup \left\{\mathsf{B}_{i,j}:  i \in \llbracket 0,\lceil {t^{-1}}\rceil \rrbracket \, ,  \, j \in \left\llbracket-\left\lceil t^{-2/3} \big(\log t^{-1}\big)^{-1/3}   \right\rceil, \left\lceil t^{-2/3} \big(\log t^{-1}\big)^{-1/3}  \right\rceil\right\rrbracket\right\}\\
&\subseteq &\bigcup\left\{\mathsf{D}_{i,j}^{(1)}\cup \mathsf{D}_{i,j}^{(2)}:i \in \llbracket 0,\lceil {t^{-1}}\rceil \rrbracket  \,   ,  \,  j \in \left\llbracket-\left\lceil t^{-2/3} \big(\log t^{-1}\big)^{-1/3}  \right\rceil, \left\lceil t^{-2/3} \big(\log t^{-1}\big)^{-1/3}  \right\rceil\right\rrbracket\right\} \, ,
\end{eqnarray*}
where~(\ref{e:imp}) was used in the latter inclusion.

Thus, with $c,k_0,n_0$ as in the statement of Theorem~ \ref{l:flucpoly}, for any fixed $t$ small enough that $\log t^{-1}\geq 2^2k_0^3$, and all $n\geq n_0(2t)^{-1}$, we have by a union bound and the translation invariance of the environment,
\begin{eqnarray*}
&& \P\left(\,  t^{-2/3} \big(\log t^{-1}\big)^{-1/3} \mathrm{MTF}_n(t)   > C \, \right)\\
&\leq & \big( 2 t^{-2/3}\big(\log t^{-1}\big)^{-1/3}+2 \big)(t^{-1}+2)\P\left(\mathrm{TF}_{n;(0,0)}^{(0,2t)}>(2^{-1}C-1)t^{2/3}\big(\log t^{-1}\big)^{1/3}\right)\\
&\leq & 2(t^{-2/3}+1)(t^{-1}+2)\exp\left\{-c(2^{-1}C-1)^3\log t^{-1}\right\}\leq  8\cdot t^{c(C/2-1)^3-5/3}.
\end{eqnarray*}
Here the second inequality follows from Theorem~\ref{l:flucpoly} with $\bm{t}=2t,\bm{k}=2^{-2/3}(2^{-1}C-1)\big(\log t^{-1}\big)^{1/3}$ and $\bm{n}=n$ being the  parameter settings. The assumptions $\log t^{-1}\geq 2^2k_0^3$, and $n\geq n_0(2t)^{-1}$ ensure that $\bm{n}\geq n_0\bm{t}^{-1}$ and $\bm{k}\geq k_0$ for any $C\geq 2$. 

Finally, choosing $C$ large enough that $c\left(C/2-1\right)^3>5/3$, we learn that
\[\P\left(\,  t^{-2/3} \big(\log t^{-1}\big)^{-1/3} \mathrm{MTF}_n(t)   > C \, \right)\rightarrow 0\quad \mbox{as } \, \, {t\searrow 0}\, , \]
whenever $n = n(t)$ verifies $n \geq n_0(2t)^{-1}$.

This completes the proof of Theorem~\ref{t:flucshortpol}.
\qed

\section{Exponent pair $(1/3,2/3)$ for polymer weight: Proof of Theorem \ref{t:holdweight}}\label{s:1/3,2/3}  
A lemma and two propositions will lead to the proof of Theorem \ref{t:holdweight} on the \holder continuity of $[1,2]\mapsto\R:t\mapsto\mathsf{Wgt}_n(t)$, the polymer weight profile under vertical displacement.

\begin{lemma}\label{l:upwtineq}
There exist positive constants $n_0,r_0,s_0,c_0$  such that, for all $n \geq n_0$, $z\in [1,2]$, $t\in [r_0n^{-1},2-z]$ and $s\in [s_0,10(nt)^{2/3}]$,
\[\P\left(|\mathsf{Wgt}_n(z+t)-\mathsf{Wgt}_n(z)|\geq st^{1/3}\right)\leq 5e^{-c_0s^{3/2}}\,.\]
\end{lemma}

We postpone the proof to Section \ref{ss:upbdwt} and first see how the lemma implies the upper bound in Theorem \ref{t:holdweight}. This bound follows from Lemma \ref{l:upwtineq} similarly to how Theorem \ref{t:main1} is derived from Proposition \ref{l:boundatz}.

\begin{proposition}\label{t:upbdwt}
The sequence $\{\mathsf{Wgt}_n\}_{n\in \N}$ is tight in $(C[1,2],\|\cdot\|_\infty)$. Moreover, if $\mathsf{Wgt}_*$ is the weak limit of a weakly converging subsequence of $\{\mathsf{Wgt}_n\}_{n\in \N}$, then there exists a positive constant~$C$ not depending on the particular weak limit $\mathsf{Wgt}_*$ such that, almost surely,
\begin{equation}\label{e:wtup}
\limsup_{t\searrow 0}\sup_{1\leq z\leq 2-t}{\left|\mathsf{Wgt}_*(z+t)-\mathsf{Wgt}_*(z)\right|} \, {t^{-1/3}\big(\log t^{-1}\big)^{-2/3}}\leq C\,.
\end{equation}
\end{proposition}

Lemma \ref{l:upwtineq} holds only for $t\in [\max\{r_0n^{-1},10^{-3/2}s^{3/2}n^{-1}\},2-z]$ for some fixed constant $r_0>0$, and not for all $t\in [0,1-z]$, as was the case in Proposition \ref{l:boundatz}. Hence, we directly show tightness in the following proof instead of applying Kolmogorov-Chentsov's tightness criterion. 


\noindent{\bf Proof of Proposition \ref{t:upbdwt}.}   To show the first statement, concerning tightness, we follow the proof of the tightness criterion used to derive \cite[Theorem~$12.3$]{Bil68}. To this end, it is enough to show  that, for given $\epsilon,\eta>0$, there exist $\delta \in [0,1]$,
which we may harmlessly suppose to verify $\delta^{-1} \in \N$, and $N_0 \in \N$ such that, for all $n\geq N_0$, 
\begin{equation}\label{e:toshowtight}
\sum_{j<\delta^{-1}} \P\left(\sup_{j\delta\leq u\leq (j+1)\delta}|\mathsf{Wgt}_n(1+u)-\mathsf{Wgt}_n(1+j\delta)|\geq \epsilon\right)<\eta\,.
\end{equation}
Assume then that $\epsilon,\eta>0$ are given small constants. For the time being, fix some $\delta>0$ small to be chosen later (depending on $\epsilon$ and $\eta$). 

Now fix any $M>1$. For any $z_1,z_2\in [1,2]$ such that $|z_1-z_2|=10^{-1}\epsilon n^{-2/3}$, set $t = \vert z_1 - z_2 \vert$. For all $\lambda\in [0,\epsilon]$, clearly $\lambda t^{-1/3}\leq 10(nt)^{2/3}$. Hence, choosing $s=\lambda t^{-1/3}$ in Lemma \ref{l:upwtineq}, one gets, for all $n$ large enough,
\begin{equation}\label{e:onepointwt}
\P\Big(|\mathsf{Wgt}_n(z_1)-\mathsf{Wgt}_n(z_2)|\geq \lambda\Big)\leq K_M\lambda^{-3M}|z_1-z_2|^M,
\end{equation}
for some constant $K_M$ depending only on $M$. 

To establish tightness, the general strategy is to bound the distribution of the maximum of certain fluctuations. To achieve this, we crucially use the bound in \eqref{e:onepointwt} together with the inequality in \cite[Theorem~$12.2$]{Bil68} that bounds the maximum of partial sums. To this end, fix
 $j<\delta^{-1}$, and break the interval $[j\delta,(j+1)\delta]$ into $\lceil \delta\beta^{-1}\rceil$-many subintervals of length $\beta:=10^{-1}\epsilon n^{-2/3}$ each, and follow the proof of the inequality in \cite[Theorem~$12.2$]{Bil68} to obtain
 \begin{equation}\label{e:wgt1}
 \P\left(\max_{0\leq i\leq \lceil \delta \beta^{-1}\rceil}|\mathsf{Wgt}_n(1+j\delta+i\beta)-\mathsf{Wgt}_n(1+j\delta)|\geq \frac{\epsilon}{2}\right)\leq K'_M\epsilon^{-3M}\delta^M,
 \end{equation}
 for some appropriate constant $K'_M$ depending only on $M$. 
 Note that by \cite[Theorem~$12.2$]{Bil68} it directly follows that if \eqref{e:onepointwt} holds for all $\lambda>0$, then \eqref{e:wgt1} holds for all $\epsilon>0$. However, in our case \eqref{e:onepointwt} holds for all $\lambda\in [0,\epsilon]$, instead of all $\lambda>0$. Hence, we resort to the proof of \cite[Theorem~$12.2$]{Bil68} which shows that if for some fixed $\epsilon>0$, \eqref{e:onepointwt} holds for all $\lambda\in [0,\epsilon]$, then \eqref{e:wgt1} holds for that particular $\epsilon$.

Now, fix any $i\in \llbracket 0,\lceil \delta \beta^{-1}\rceil-1\rrbracket$. 
For any $u\in [j\delta+i\beta,j\delta+(i+1)\beta]$, it clearly follows from the definition~\eqref{e:defweight},
 \[\weight_{n;(0,1+u)}^{(0,1+j\delta+(i+1)\beta)}\geq -2n^{2/3}(1+j\delta+(i+1)\beta-(1+u))\geq -2n^{2/3}\beta\,,\,\,\mbox{and}\]
\[\weight_{n;(0,1+j\delta+i\beta)}^{(0,1+u)}\geq -2n^{2/3}(1+u-(1+j\delta+i\beta))\geq -2n^{2/3}\beta\,.\]
Thus, for any $u\in [j\delta+i\beta,j\delta+(i+1)\beta]$, by superaddivity of polymer weights described in \eqref{e:superadd},
\[
\weight_{n;(0,0)}^{(0,1+j\delta+i\beta)}-2n^{2/3}\beta\leq \weight_{n;(0,0)}^{(0,1+j\delta+i\beta)}+\weight_{n;(0,1+j\delta+i\beta)}^{(0,1+u)}\leq \weight_{n;(0,0)}^{(0,1+u)}\,\,\, \mbox{and}\]
\[
\weight_{n;(0,0)}^{(0,1+u)}\leq \weight_{n;(0,0)}^{(0,1+j\delta+(i+1)\beta)}-\weight_{n;(0,1+u)}^{(0,1+j\delta+(i+1)\beta)}\leq \weight_{n;(0,0)}^{(0,1+j\delta+(i+1)\beta)}+2n^{2/3}\beta\,.
\]
This, together with \eqref{e:Wmod}, imply that for any $i\in \llbracket 0,\lceil \delta \beta^{-1}\rceil-1\rrbracket$ and $u\in [j\delta+i\beta,j\delta+(i+1)\beta]$,
\begin{eqnarray}\label{e:wgt2}
n^{1/3}\left|\mathsf{Wgt}_n(1+u)-\mathsf{Wgt}_n(1+j\delta)\right|\leq  2n\beta+2+ n^{1/3}\max\Big\{\left|\mathsf{Wgt}_n(1+j\delta+i\beta)-\mathsf{Wgt}_n(1+j\delta)\right|,\nonumber\\
  \left|\mathsf{Wgt}_n(1+j\delta+(i+1)\beta)-\mathsf{Wgt}_n(1+j\delta)\right|\Big\}.
\end{eqnarray}

Since $2n\beta= 5^{-1}\epsilon n^{1/3}$, for all $n$ large enough that $2n^{-1/3}\leq \epsilon/5$, \eqref{e:wgt1} and \eqref{e:wgt2} imply 
\begin{eqnarray*}
&&\P\left(\sup_{j\delta\leq u\leq (j+1)\delta}\left|\mathsf{Wgt}_n(1+u)-\mathsf{Wgt}_n(1+j\delta)\right|\geq \epsilon\right)\\
&\leq& \P\left(\max_{0\leq i\leq \lceil \delta \beta^{-1}\rceil}|\mathsf{Wgt}_n(1+j\delta+i\beta)-\mathsf{Wgt}_n(1+j\delta)|\geq \frac{\epsilon}{2}\right)\leq K'_M\epsilon^{-3M}\delta^M.
\end{eqnarray*}
Thus, by choosing $\delta$ small enough that $K'_M\epsilon^{-3M}\delta^{M-1}<\eta$, we obtain \eqref{e:toshowtight}, and hence tightness.

To show \eqref{e:wtup}, we follow the proof of Theorem \ref{t:main1}(b). 
 Let $n_0,r_0,s_0$ and $c_0$ be as in Lemma~\ref{l:upwtineq}. For any fixed $m\in \N$ such that $c_1\left(\log 2^m\right)^{2/3}\geq s_0$, and any $j\in \{0,1,2,\cdots,2^m-1\}$, and all $n\geq \max\{r_02^m,10^{-3/2}c_1^{3/2}2^m\log 2^m\}$, by applying Lemma~\ref{l:upwtineq} with the parameters $\bm{n}=n,\bm{t}=2^{-m}$ and $\bm{s}=c_1\left(\log 2^m\right)^{2/3}$, it follows that 
\[\P\Bigg( \big\vert \mathsf{Wgt}_n(1+(j+1)2^{-m})-\mathsf{Wgt}_n(1+j2^{-m}) \big\vert > c_12^{-\frac{m}{3}}\left(\log 2^m\right)^{2/3}\Bigg) \, \leq \, 5\cdot 2^{-m(c_0c_1^{3/2})}\,.\]
Now, observe that \eqref{e:main1} in the proof of Lemma~\ref{l:dyad} carries over verbatim to the present case. By choosing $c_1$ high enough  that $c_0c_1^{3/2}>1$, and exactly imitating the rest of the proof of Lemma \ref{l:dyad} followed by the proof of Theorem \ref{t:main1}(b), we complete the proof of Proposition~\ref{t:upbdwt}.
\qed

Turning to prove the lower bound in \eqref{e:wgtholdbounds}, we restate it now.
\begin{proposition}\label{p:wgtlwb} There exists a constant $c>0$ such that, almost surely,
\[\liminf_{t\searrow 0}\sup_{1\leq z\leq 2-t} \, t^{-1/3}\big(\log t^{-1}\big)^{-2/3}  \big\vert \mathsf{Wgt}_*(z+t)-\mathsf{Wgt}_*(z) \big\vert \,  \geq \,  c\,.\]
\end{proposition}
 This result will follow directly from weight superadditivity, i.e. $\weight_{n;(0,0)}^{(0,1+z+t)}-\weight_{n;(0,0)}^{(0,1+z)}\geq \weight_{n;(0,1+z)}^{(0,1+z+t)}$ for $z,t>0$,
 control on weight with given endpoints via Theorem~\ref{t:moddevlow},
  independence in disjoint strips, and the weight $\weight_{n;(0,1+z)}^{(0,1+z+t)}$ depending on the configuration in the strip delimited by the lines $y=1+z$ and $y=1+z+t$.  The proof is reminiscent of an argument for a similar statement made for Brownian motion: see the proof on page $362$ of Exercise $1.7$ in the book \cite{MP}. 


\noindent{\bf Proof of Proposition \ref{p:wgtlwb}.} We need to show that, for some constant $c>0$, almost surely, there exists $\epsilon>0$ such that, for all $0<t<\epsilon$ and some $z\in[1,2-t]$,
\[|\mathsf{Wgt}_*(z+t)-\mathsf{Wgt}_*(z)|\geq ct^{1/3}\big(\log t^{-1}\big)^{2/3}.\]

Let $c>0$ satisfy $2^{3/2}c_2c^{3/2}<1$, where $c_2$ arises from  Theorem \ref{t:moddevlow}. For integers $n,m\geq 1$ and $k\in \{0,1,2,\cdots,m-1\}$, we define the events
\[\mathsf{A}_{k,m,n}=\left\{\mathsf{Wgt}_n\left(1+{(k+1)}{m}^{-1}\right)-\mathsf{Wgt}_n\left(1+{k}{m}^{-1}\right)\geq c{m^{-1/3}}\left(\log m\right)^{2/3}\right\} \]
and
\[\mathsf{A}_{k,m}=\left\{\mathsf{Wgt}_*\left(1+{(k+1)}{m}^{-1}\right)-\mathsf{Wgt}_*\left(1+{k}{m}^{-1}\right)\geq c{m^{-1/3}}\left(\log m\right)^{2/3}\right\}\,.\]
Also let 
\[\mathsf{B}_{k,m,n}=\left\{\weight_{n;(0,1+km^{-1})}^{(0,1+(k+1)m^{-1})}\geq c{m^{-1/3}}\left(\log m\right)^{2/3}+2n^{-1/3}\right\}\,.\]

Let $n_0,s_0$ and $c_2$ be as in Theorem \ref{t:moddevlow}, and let $m_0$ be large enough that $2c(\log m_0)^{2/3}\geq \max\{s_0,4n_0^{-1/3}\}$. Then from Theorem \ref{t:moddevlow} with parameter settings $\bm{t_1}= 1+km^{-1},\bm{t_2}=1+(k+1)m^{-1},\bm{\tot}=m^{-1},\bm{n}=n$ and $\bm{s}=2c(\log m)^{2/3}$, for all $m\geq m_0$ and $n\geq n_0m$, 
\begin{equation}\label{e:Bomn}
\P(\mathsf{B}_{0,m,n})\geq \P\left(\weight_{n;(0,1+km^{-1})}^{(0,1+(k+1)m^{-1})}\geq 2c{m^{-1/3}}\left(\log m\right)^{2/3}\right)\geq  e^{-2^{3/2}c_2c^{3/2}\log m}=m^{-2^{3/2}c_2c^{3/2}}\,.
\end{equation}
Here the first inequality follows because
\[cm^{-1/3}\left(\log m\right)^{2/3}\geq cm^{-1/3}\left(\log m_0\right)^{2/3}\geq 2n_0^{-1/3}m^{-1/3}\geq 2n^{-1/3}\,\]
for $m\geq m_0$, $2c(\log m_0)^{2/3}\geq 4n_0^{-1/3}$ and $n\geq n_0m$.

Now $\mathsf{B}_{k,m,n}$ are i.i.d. random variables for $k\in\{0,1,2,\cdots,m-1\}$ as the weights of polymers over disjoint regions are independent. Also using  $\weight_{n;(0,0)}^{(0,1+(k+1)m^{-1})}-\weight_{n;(0,0)}^{(0,1+km^{-1})}\geq \weight_{n;(0,1+km^{-1})}^{(0,1+(k+1)m^{-1})}$ by superadditivity of polymer weights, together with \eqref{e:Wmod}, we get that $\mathsf{B}_{k,m,m}\subseteq \mathsf{A}_{k,m,n}$.  Thus, using \eqref{e:Bomn}, for all $m\geq m_0$ and $n\geq n_0m$, 
\begin{eqnarray}\label{e:Acom}
\P\left(\bigcap_{k=0}^{m-1} \mathsf{A}^c_{k,m,n}\right)&\leq& \P\left(\bigcap_{k=0}^{m-1} \mathsf{B}^c_{k,m,n}\right)
=(1-\P(\mathsf{B}_{0,m,n}))^m\nonumber\\
&\leq & \exp\left\{-m\P(\mathsf{B}_{0,m,n})\right\}\leq \exp\Big\{-m^{1-2^{3/2}c_2c^{3/2}}\Big\}\,,
\end{eqnarray}
where we use that $1-x\leq e^{-x}$ for all $x\geq 0$.


Next, similarly to the first part of the proof of Lemma \ref{l:dyad}, let $\{\mathsf{Wgt}_{n_r}\}_r$ be a subsequence of $\{\mathsf{Wgt}_{n}\}_n$ such that $\mathsf{Wgt}_{n_r}\Rightarrow \mathsf{Wgt}_*$ as random variables in $(C[1,2],\|\cdot\|_\infty)$ (where $\Rightarrow$ denotes convergence in distribution).
 Since for $a,b\in [1,2]$, the map $T_{a,b}$ defined by $(C[1,2],\|\cdot\|_\infty)\mapsto (\R,|\cdot|):f\mapsto f(a)-f(b)$ is continuous, the set \[\mathsf{U}:=\bigcap\left\{T^{-1}_{1+(k+1)m^{-1},1+km^{-1}}\left(-\infty,cm^{-1/3}\big(\log m\big)^{2/3}\right):k=0,1,\cdots,m-1\right\}\]
 is open. Thus, by the Portmanteau theorem,
 \[\P\left(\bigcap_{k=0}^{m-1} \mathsf{A}^c_{k,m}\right)\leq \liminf_r\P\left(\bigcap_{k=0}^{m-1} \mathsf{A}^c_{k,m,n_r}\right)\leq \limsup_n\P\left(\bigcap_{k=0}^{m-1} \mathsf{A}^c_{k,m,n}\right)\,.\]
 
From here, using \eqref{e:Acom} and that our given choice of the constant $c$ ensures $2^{3/2}c_2c^{3/2}<1$, we get
\begin{eqnarray*}
\sum_{m=m_0}^\infty \P\left(\bigcap_{k=0}^{m-1} \mathsf{A}^c_{k,m}\right)\leq \sum_{m=m_0}^\infty \limsup_n\P\left(\bigcap_{k=0}^{m-1} \mathsf{A}^c_{k,m,n}\right)\leq \sum_{m=m_0}^\infty \exp\Big\{-m^{1-2^{3/2}c_2c^{3/2}}\Big\} <\infty\,.
\end{eqnarray*}
Hence, using the Borel-Cantelli lemma, almost surely there exists $M_0\in \N$ such that for all $m\geq M_0$, one has some $k_m\leq m-1$ with $z=1+{k_m}{m}^{-1}$ satisfying
\[\left|\mathsf{Wgt}_*(z+{m}^{-1})-\mathsf{Wgt}_*(z)\right|\geq c{m^{-1/3}}\left(\log m\right)^{2/3}.\]
Let $\epsilon= M_0^{-1}$. Also let $M_0^{-1}$ be small enough in the sense of Proposition \ref{t:upbdwt}: namely, almost surely for all $t\in [0,M_0^{-1}]$, $\sup_{1\leq z\leq 2-t}{\left|\mathsf{Wgt}_*(z+t)-\mathsf{Wgt}_*(z)\right|}{t^{-1/3}\big(\log t^{-1}\big)^{-2/3}}\leq 2C$. Then, for any given $t\in [0,\epsilon]$, let $m$ be such that ${(m+1)^{-1}}<t\leq {m}^{-1}$. Then for $z=1+{k_m}{m}^{-1}$, 
\begin{eqnarray*}
&&|\mathsf{Wgt}_*(z+t)-\mathsf{Wgt}_*(z)|\\
&\geq &\left|\mathsf{Wgt}_*\left(z+{m}^{-1}\right)-\mathsf{Wgt}_*(z)\right|-\left|\mathsf{Wgt}_*(z+t)-\mathsf{Wgt}_*\left(z+{m}^{-1}\right)\right|\\
&\geq&c{m^{-1/3}}\left(\log m\right)^{2/3}-2C\left({m}^{-1}-{(m+1)}^{-1}\right)^{1/3}\Big(\log\left({m}^{-1}-{(m+1)}^{-1}\right)^{-1}\Big)^{2/3}.
\end{eqnarray*}
As the second term decays much faster than the first, choosing $M_0$ large enough so that the second term is smaller that $2^{-1} c {m^{-1/3}}(\log m)^{2/3}$ gives the result.
\qed

\noindent{\bf Proof of Theorem \ref{t:holdweight}.}
This result follows from Proposition \ref{t:upbdwt} and Proposition \ref{p:wgtlwb}.
\qed

\subsection{Upper bound on polymer weight fluctuation: Proof of Lemma \ref{l:upwtineq}}\label{ss:upbdwt}

In this subsection, we complete the proof of Theorem \ref{t:holdweight}. 
The remaining element, Lemma \ref{l:upwtineq}, will be derived from  Lemmas~\ref{l:upineq1} and~\ref{l:upwtineq2}.




\begin{lemma}\label{l:upineq1} There exist  positive constants $s_0,r_0$ and $c_0$ such that for $s\geq s_0$, $z\in [1,2]$ and $t\in [r_0n^{-1},2-z]$, 
\[\P\left(\weight_{n;(0,0)}^{(0,z)}\geq \weight_{n;(0,0)}^{(0,z+t)}+st^{1/3}\right)\leq e^{-c_0s^{3/2}}.\]
\end{lemma}
\noindent{\bf Proof.} Using $\weight_{n;(0,0)}^{(0,z+t)}\geq \weight_{n;(0,0)}^{(0,z)}+\weight_{n;(0,z)}^{(0,z+t)}$,
we see that, for $nt\geq r_0$ and $s\geq s_0$,
\[\P\left(\weight_{n;(0,0)}^{(0,z)}\geq \weight_{n;(0,0)}^{(0,z+t)}+st^{1/3}\right)\leq \P\big(\weight_{n;(0,z)}^{(0,z+t)}\leq -st^{1/3}\big)\leq e^{-cs^{3/2}},\]
where the latter inequality follows from the moderate deviation estimate Theorem~\ref{t:moddev}, with $\bm{t_1}=z,\bm{t_2}=z+t,\bm{n}=n$ and $\bm{s}=s$, and setting $r_0$ and $s_0$ to equal $n_0$ and $s_0$ respectively from the statement of Theorem~\ref{t:moddev}. 
\qed

\begin{figure}[h]
\centering
\includegraphics[width=0.9\textwidth]{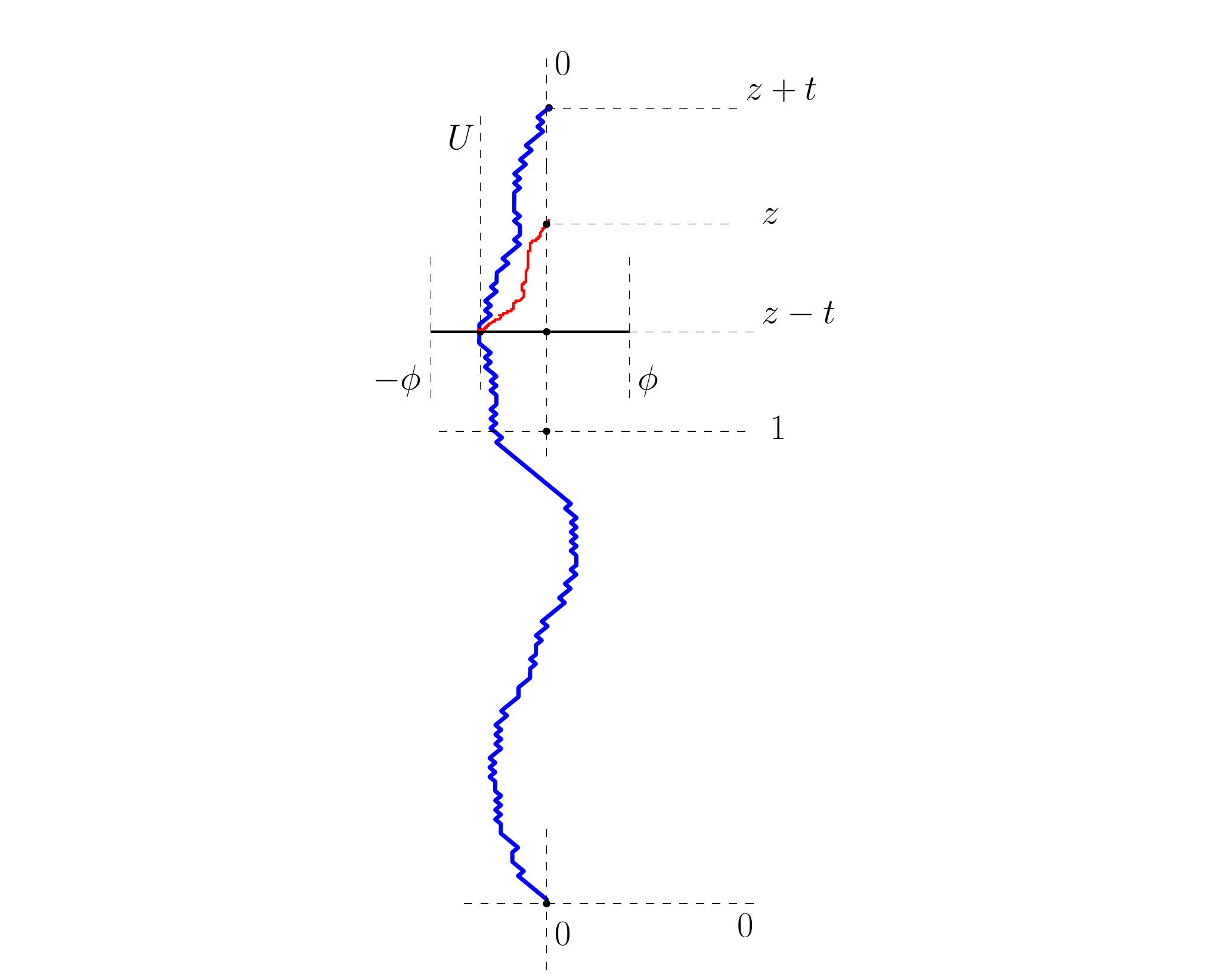}
\caption{When the
thick blue polymer $\rho_{n;(0,0)}^{\leftarrow;(z+t,0)}$ crosses height $z-t$ without immoderately high fluctuation, it may be diverted via the red polymer to form a path of comparable weight from $(0,0)$ to $(z,0)$.}
\label{f:fig2weight}
\end{figure}

Next is the more subtle of the two constituents of Lemma \ref{l:upwtineq}.
\begin{lemma}\label{l:upwtineq2}There exist  positive constants $n_0,s_2,r_1$ and $c_0$ such that, for $n\geq n_0$, $t\in [r_1n^{-1},2-z]$, $s\in [s_2,10(nt)^{2/3}]$ and $z\in [1,2]$,
\begin{equation}\label{e:wgtuplem}
\P\left(\weight_{n;(0,0)}^{(0,z+t)}\geq \weight_{n;(0,0)}^{(0,z)}+st^{1/3}\right)\leq 4 \, e^{-c_0s^{3/2}}.
\end{equation}
\end{lemma}
This proof is reminiscent of arguments used in \cite{BSS14} and \cite{BSS17++}. We first explain the basic idea, which is illustrated in Figure~\ref{f:fig2weight}.
A path may be formed from $(0,0)$
to $(0,z)$ by following the route of a polymer from $(0,0)$ to $(0,z+t)$ until its location, $(U,z-t)$ say, at height $z-t$; and then following a polymer from $(U,z-t)$ to $(0,z)$. The discrepancy in weight between the original polymer, from $(0,0)$ to $(0,z+t)$, and the newly formed path, from $(0,0)$ to $(0,z)$, is equal to the difference in weights between the polymer from $(U,z-t)$ to $(0,z+t)$
and that from $(U,z-t)$ to $(0,z)$. The latter two polymers have duration of order $t$; Theorem~\ref{t:unifmoddev} may then show that  their weights have order $t^{1/3}$. Thus, the weight difference $\weight_{n;(0,0)}^{(0,z+t)} - \weight_{n;(0,0)}^{(0,z)}$, which is at most the discrepancy we are considering, is seen to be unlikely to exceed order $t^{1/3}$.


\noindent{\bf Proof of Lemma \ref{l:upwtineq2}.}
To implement this idea, we will consider, for definiteness, the leftmost polymer from $(0,0)$ to $(z+t,0)$, namely $\rho_{n;(0,0)}^{\leftarrow;(z+t,0)}$.
 In accordance with the notation in the plan, we will set
 $U = \rho_{n;(0,0)}^{\leftarrow;(z+t,0)}(z-t)$.
 
The height-$(z-t)$ polymer location $U$ typically has order $t^{2/3}$. The plan will run into trouble if $U$ is atypically high, because then the two short polymers running to $(0,z+t)$ and $(0,z)$ from $(U,z-t)$
will have large negative weights dictated by parabolic curvature. 

To cope with this difficulty, we introduce a {\em good} event $\mathsf{G}$,
\[\mathsf{G}=\left\{\left|U\right|\leq \phi\right\}\,,\]
specified in terms of a parameter 
 $\phi$ that is set equal to $D^{-1}s^{1/2}(2t)^{2/3}$. Here,  the  constant $D$ is chosen to be $2^{2/3}10^{1/2}{C_0}$, with  $C_0$ given by Theorem~\ref{t:unifmoddev}. 
In view of Theorem~\ref{t:carsestimate2}, this choice of $\phi$
ensures that the event $\mathsf{G}$ fails to occur with probability of order $\exp \big\{ - \Theta(1) s^{3/2} \big\}$. (The appearance of the factor of $D^{-1}$ in $\phi$ is a detail concerning values of $s$ in Lemma~\ref{l:upwtineq2} close to the maximum value $10 (nt)^{2/3}$. ) 
 

Indeed, applying Theorem \ref{t:carsestimate2} with $\bm{n}=n,\bm{t_1}=0,\bm{t_2}=z+t,\bm{t}=z-t,\bm{x}=0,\bm{y}=0$ and $\bm{s}=D^{-1}s^{1/2}$, we find that, when $n\geq n_0$ (a bound which ensures that the hypothesis that $\bm{n\tot}
\geq n_0$ is met) and $s \geq s_1$, 
\begin{equation}\label{e:boundGcom}
\P(\mathsf{G}^c) \, \leq \,  2 \, \exp \big\{ -c D^{-3} s^{3/2} \big\} \, ,
\end{equation}
where the positive constants $c$ and $s_1$
are provided by the theorem being applied.

When $\mathsf{G}$ occurs, 
\[|U|\leq D^{-1}s^{1/2}(2t)^{2/3}\leq D^{-1}2^{2/3}10^{1/2}tn^{1/3}<tn^{1/3}\,,\]
because $s\leq 10(nt)^{2/3}$, $D = 2^{2/3} 10^{1/2} C_0$ and $C_0>1$. As we saw in Subsection~\ref{sss:comp}, it is this bound on $\vert U \vert$
that ensures the existence of polymers between $(U,z-t)$ and $(0,z)$. By superadditivity of polymer weights, we thus have
\[\weight_{n;(0,0)}^{(0,z)}\geq \weight_{n;(0,0)}^{(U,z-t)}+\weight_{n;(U,z-t)}^{(0,z)}\,.\]

 Thus, when $\mathsf{G}$ occurs,
 \begin{eqnarray*}
 \weight_{n;(0,0)}^{(0,z+t)}- \weight_{n;(0,0)}^{(0,z)}&\leq& \weight_{n;(0,0)}^{(0,z+t)}-\weight_{n;(0,0)}^{(U,z-t)}-\weight_{n;(U,z-t)}^{(0,z)}\\
 &=&\weight_{n;(U,z-t)}^{(0,z+t)}-\weight_{n;(U,z-t)}^{(0,z)}
 \, \leq \, \sup_{x\in [-\phi,\phi]}\left(\weight_{n;(x,z-t)}^{(0,z+t)}-\weight_{n;(x,z-t)}^{(0,z)}\right)\,,
 \end{eqnarray*}
 where the equality is dependent on the definition of $U$ and the final inequality on the occurrence of $\mathsf{G}$.
 We see then that 
\begin{eqnarray}
&&\P\left(\mathsf{G}\cap \left\{\weight_{n;(0,0)}^{(0,z+t)}\geq \weight_{n;(0,0)}^{(0,z)}+st^{1/3}\right\}\right) \nonumber \\
&\leq&  \P\left(\sup_{x\in [-\phi,\phi]}\left(\weight_{n;(x,z-t)}^{(0,z+t)}-\weight_{n;(x,z-t)}^{(0,z)}\right)\geq st^{1/3}\right)\nonumber\\
& \leq & \P\left(\sup_{x\in [-\phi,\phi]}\left|\weight_{n;(x,z-t)}^{(0,z+t)}\right|>2^{-1}st^{1/3}\right)+ \P\left(\sup_{x\in [-\phi,\phi]}\left|\weight_{n;(x,z-t)}^{(0,z)}\right|>2^{-1}st^{1/3}\right) \label{e:ineqgcap} \,.
\end{eqnarray}

The latter two probabilities will each be bounded above by a union bound over several applications of Theorem~\ref{t:unifmoddev}. Addressing the first of these probabilities to begin with, we set parameters for a given application of  the theorem, taking
$\bm{I}$ to be a given interval of length at most $t^{2/3}$
contained in $[-\phi,\phi]$ and $\bm{J} = \{ 0 \}$,
and also setting
 $\bm{n}=n, \bm{t_1}=z-t,\bm{t_2}=z$ and $\bm{s}=4^{-1}s$.

The theorem's hypothesis concerning inclusion for the interval ${\bf I}$
(and ${\bf J}$) is ensured because  
\[|x|\leq D^{-1}s^{1/2}(2t)^{2/3}\leq 2^{2/3}10^{1/4}D^{-1}s^{1/4}n^{1/6}t^{5/6}< C_0^{-1}s^{1/4}n^{1/6}t^{5/6}\,,\]
for $x\in [-\phi,\phi]$, where here we use $s\leq 10(nt)^{2/3}$ and $D= 2^{2/3}10^{1/2}C_0>2^{2/3}10^{1/4}C_0$. 

In these applications of Theorem~\ref{t:unifmoddev},
the parabolic curvature term inside the supremum,  $t^{-4/3} x^2$, is at most $t^{-4/3} \phi^2$.
It is thus also at most $s/4$, because $\phi = D^{-1} s^{1/2}(2t)^{2/3}$ and $D \geq 2^{5/3}$.

Thus, dividing $[-\phi,\phi]$ into $\lceil 2^{5/3}D^{-1}s^{1/2}\rceil$-many consecutive intervals of length at most $t^{2/3}$, we are indeed able to apply Theorem \ref{t:unifmoddev} and a union bound, finding that, for $n_0 \in \N$ and $C,c  > 0$ the constants furnished by the theorem, and for $nt \geq n_0$,
\begin{eqnarray*}
&&\P\left(\sup_{x\in [-\phi,\phi]}\left|\weight_{n;(x,z-t)}^{(0,z)}\right|>2^{-1}st^{1/3}\right)\\
&\leq &\P\left(\sup_{x\in [-\phi,\phi]}\left|t^{-1/3}\weight_{n;(x,z-t)}^{(0,z)}+t^{-4/3}x^2\right|>4^{-1}s\right)\leq \lceil 2^{5/3}D^{-1}s^{1/2}\rceil C e^{-cs^{3/2}}\leq e^{-c's^{3/2}}\,,
\end{eqnarray*}
for $c'=2^{-1}c$ and $s\geq s_0$ where $s_0$ is chosen in such a way that $e^{2^{-1}cs_0^{3/2}}\geq C \lceil 2^{5/3}D^{-1}s_0^{1/2}\rceil$. 

The second probability in \eqref{e:ineqgcap} is bounded above by similar means. Several applications of Theorem~\ref{t:unifmoddev} will be made. In a given application, the parameters $\bm{I}, \bm{J}, \bm{n}$ and $\bm{s}$ are chosen as before, but we now set $\bm{t_1}=z-t$ and $\bm{t_2}=z+t$, so that $\bm{\tot}$ equals $2t$, rather than $t$. 
The curvature term $(2t)^{-4/3} x^2$ is bounded above by  $(2t)^{-4/3} \phi^2$, a smaller bound than before, so that the preceding bound of  $s/4$ remains valid. The condition for inclusion for  the intervals  $\bm{I}$ (and $\bm{J}$), namely  $\phi \leq C_0^{-1}s^{1/4}n^{1/6}(2t)^{5/6}$, is weaker than it was previously and is thus satisfied. Hence, using Theorem \ref{t:unifmoddev} and a union bound, we find that, for all $n\geq 2^{-1}n_0 t^{-1}$,
\[\P\left(\sup_{x\in [-\phi,\phi]}\left|\weight_{n;(x,z-t)}^{(0,z+t)}\right|>2^{-1}st^{1/3}\right)\leq e^{-c's^{3/2}}\,,\]
for $s\geq s_0$.

Combining \eqref{e:boundGcom} and \eqref{e:ineqgcap} with the two bounds just derived, we obtain Lemma~\ref{l:upwtineq2} by taking $c_0 > 0$
to be less than $\min\{c D^{-3},c'\}$, $s_2$ to be suitably greater than $\max\{s_0, s_1\}$, and $r_1=2^{-1}n_0$. 
\qed

\noindent{\bf Proof of Lemma \ref{l:upwtineq}.}  This follows immediately using \eqref{e:Wmod} and from 
Lemmas~\ref{l:upineq1} and~\ref{l:upwtineq2} and a union bound.
\qed

\section{Lower bound on transversal fluctuation: Proof of Proposition \ref{p:lowbndTF}}\label{s:last}
In this last section we shall prove the lower bound on the transversal fluctuation of the polymer, the corresponding upper bound of which was proved in \cite[Theorem $11.1$]{BSS14} (and is stated here, with the optimal exponent in the bound, as Theorem \ref{t:transversal}). 
In fact, Proposition \ref{p:lowbndTF} does slightly more than just providing a corresponding lower bound on the quantity whose upper bound is proved in Theorem \ref{t:transversal}. Indeed, in Proposition \ref{p:lowbndTF}, one takes the minimum over the transversal fluctuations of all the polymers between two fixed points, and not just the transversal fluctuation of the leftmost one. The proof of Proposition \ref{p:lowbndTF} crucially uses the polymer weight lower tail Theorem~\ref{t:moddevlow}. We also fix the constant $\alpha_0$ in this Proposition \ref{p:lowbndTF} as $\alpha_0=C_0^{-2}3^{-5/3}10^{-1/2}$, where $C_0$ is as in Theorem~\ref{t:unifmoddev}. This choice of $\alpha_0$ ensures that the condition in the hypothesis of Theorem~\ref{t:unifmoddev} is met whenever it is applied.



\noindent{\bf Proof of Proposition \ref{p:lowbndTF}.} We prove the proposition for $t_1=0$ and $t_2=1$. The case for general $t_1<t_2$ follows readily using the scaling principle \eqref{e:screl}.

A box is a subset of $\R^2$ of the form $[a,b] \times [r_1,r_2]$, where $a \leq b$ and $r_1 \leq r_2$.
Any box has a lower and an upper side, namely $[a,b] \times \{ r_1 \}$ and $[a,b] \times \{ r_2 \}$.

The key box for the proof is $\strip$, now specified to be $[-s,s] \times [0,1]$. Proposition~\ref{p:lowbndTF} is, after all, a lower bound on the probability that there exists a polymer between $(0,0)$ and $(0,1)$ that escapes $\strip$.

\begin{figure}
\centering
\includegraphics[scale=.6]{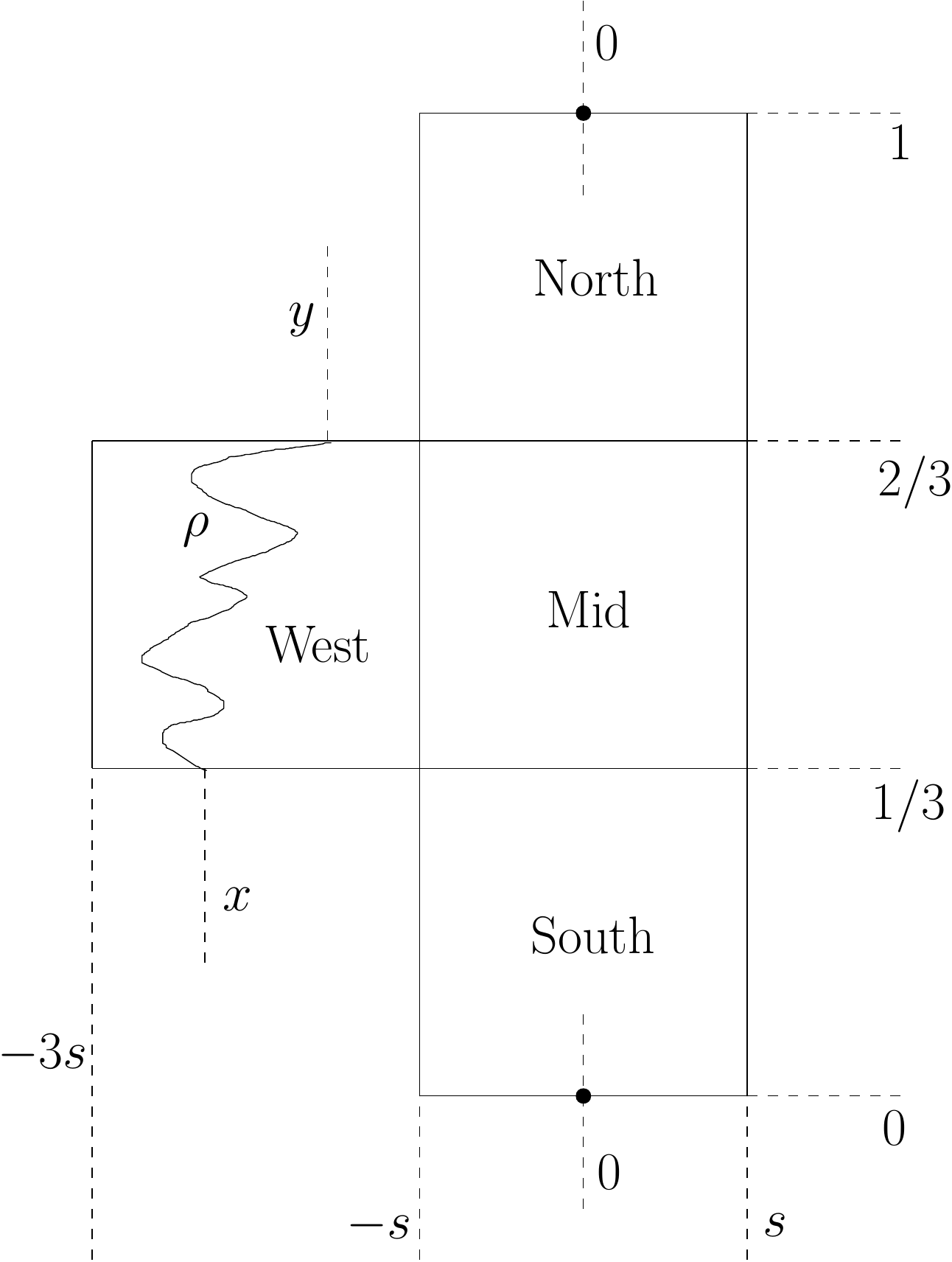} 
\caption{In Case High, the high weight path $\rho$
is extended to form a path from $(0,0)$ to $(0,1)$
whose weight exceeds that of any path between these points that remains in $\strip = \north \cup \midd \cup \south$.}\label{f.casehigh}
\end{figure}

We divide $\strip$ into three further boxes, writing
$\midd$ for the box $[-s,s] \times [1/3,2/3]$, and 
 $\south$ and  $\north$ for the boxes obtained from $\midd$ by vertical translations of $-1/3$ and $1/3$. We further set $\west$ to be the box obtained from $\midd$ by a horizontal translation of $-2s$. See Figure \ref{f.casehigh}.
 
Recall that, when $(x,t_1)$ and $(y,t_2)$ verify $n^{1/3} \tot \geq |y - x|$, we denote the polymer weight with this pair of endpoints by $\weight_{n;(x,t_1)}^{(y,t_2)}$. We now use a set theoretic notational convention to refer in similar terms to the set of weights of polymers between two collections of endpoint locations. Indeed, let $I$ and $J$ be compact real intervals. We will write
$$
\weight_{n;(I,t_1)}^{(J,t_2)} = \Big\{ 
\weight_{n;(x,t_1)}^{(y,t_2)} : x \in I, y \in J \Big\} \, ;
$$
we will ensure that whenever this notation is used, $(x,t_1) \overset{n}{\preceq} (y,t_2)$ for all $x \in I$ and $y \in J$ in the sense of Subsection \ref{sss:comp}. 
 When an interval is a singleton, $I = \{ x \}$ say, we write $(x,t_1)$ instead of $(\{x\},t_1)$ when using this notation. 

To any box $B$ and $s \in \R$,
we define the event $\high(B,s)$
that the weight of some path that
is contained in $B$ with starting point in the lower side of $B$ and ending point in the upper side of $B$ is at least $s$. 

Our approach to proving Proposition~\ref{p:lowbndTF}
gives a central role to the event $\high(\midd,300s^2)$.
It may be expected that the order of probability of this event is $\exp \big\{ - \Theta(1) s^3  \big\}$, but we do not attempt to prove this. Rather, we analyse two cases, called {\em High} and {\em Low}, according to the value of the event's probability.

We will quantify the notion of high or low probability for $\high(\midd,300s^2)$ in terms of the decay rate for a very high weight polymer between $(0,0)$ and $(0,1)$.
Indeed, noting from Theorem \ref{t:moddevlow} that there exists $C > 0$ 
such that, for $s\geq s_0$,
\begin{equation}\label{e.wthou}
 \P \Big( \weight_{n;(0,0)}^{(0,1)} \geq 1000s^2 \Big) \geq \exp \big\{ - C s^3 \big\} \, ,
\end{equation}
we declare that Case High occurs if 
$$
 \P \big( \high(\midd,300s^2) \big) \geq \exp \big\{ -2C s^3 \big\} \, ;
$$

Case Low occurs when Case High does not.

In order to analyse Case High, we introduce a {\em favourable} event $\mathsf{F}$. The event is specified as the intersection of the following events:
\begin{itemize}
\item  $\mathsf{G}_1=\left\{\inf\weight_{n;(0,0)}^{([-3s,-s],1/3)} \geq - 50s^2\right\}$;
\item $\mathsf{G}_2=\left\{\inf\weight_{n;([-3s,-s],2/3)}^{(0,1)} \geq - 50s^2\right\}$;
\item $\mathsf{G}_3=\left\{\sup \weight_{n;(0,0)}^{([-s,s],1/3)} \leq  50s^2\right\}$;
\item $\mathsf{G}_4=\left\{\sup \weight_{n;([-s,s],2/3)}^{(0,1)} \leq  50s^2\right\}$;
\item and $\mathsf{G}_5$ is the event that $\high(\midd,50s^2)$ does not occur.
\end{itemize}
Thus, the occurrence of $\mathsf{F}$ forces the absence of any high weight path inside $\midd$ that crosses this box from its lower to its upper side, while also ensuring that any polymer connecting $(0,0)$ (or $(0,1)$)
to the lower (or upper) sides of $\midd$ and $\west$ is not of very low weight. We claim that $\mathsf{F}$ is a high probability event, proving this by applying Theorem~\ref{t:unifmoddev}. Indeed, for the events $\mathsf{G}_1$ and $\mathsf{G}_3$ entailed by $\mathsf{F}$, we make several applications of Theorem~ \ref{t:unifmoddev}. For a given application, we consider the parameter settings $\bm{n}=n, \bm{t_1}=0,\bm{t_2}=1/3,
 \bm{s}=10s^2, \bm{I}=\{0\}$ and 
 $$
 \bm{J}= \big[-3s+(i-1)3^{-2/3},\max\{-3s+i3^{-2/3},s\} \big]
 $$ for some $i\in \{1,2,\cdots,\lceil4\cdot 3^{2/3}s\rceil\}$. The condition on inclusion for the intervals $\bm{I}$ and $\bm{J}$ is satisfied since for $y\in [-3s,s]$,  
\[|y|\leq 3s\leq \bm{s}^{1/2}\leq 10^{1/4}\alpha_0^{1/2}\bm{n}^{1/6}\bm{s}^{1/4}\leq 3^{5/6}10^{1/4}\alpha_0^{1/2}\bm{n}^{1/6}\bm{s}^{1/4}\bm{\tot}^{5/6}=C_0^{-1}\bm{n}^{1/6}\bm{s}^{1/4}\bm{\tot}^{5/6}\,,\]
where we use that $s\leq \alpha_0n^{1/3}$ and our given choice of $\alpha_0$ has been made so that $\alpha_0=C_0^{-2}3^{-5/3}10^{-1/2}$. Also the parabolic curvature inside the supremum is
\[\sup_{y\in [-3s,s]}3^{4/3}y^2\leq 3^{4/3}\cdot 3^2s^2<40s^2\,.\]
Thus, dividing $[-3s,s]$ into $\lceil4\cdot 3^{2/3}s\rceil$-many intervals of length at most $3^{-2/3}$ and using Theorem~\ref{t:unifmoddev} and a union bound, it follows that, for $s$ large enough and $n\geq 3n_0$,
\[\P(\mathsf{G}_1^c\cup \mathsf{G}_3^c)\leq \P\left(\sup_{y\in [-3s,s]}\left|3^{1/3}\weight_{n;(0,0)}^{(y,1/3)}+3^{4/3}y^2\right|>10s^2\right)\leq \lceil 4\cdot 3^{2/3}s\rceil Ce^{-cs^3}\leq 6^{-1}\,.\]
Similarly for the events $\mathsf{G}_2$ and $\mathsf{G}_4$, in a given application of Theorem~\ref{t:unifmoddev}, we set the parameters $\bm{n}=n,\bm{t_1}=2/3,\bm{t_2}=1,
 \bm{s}=10s^2,\bm{I}=[-3s+(i-1)3^{-2/3},\max\{-3s+i3^{-2/3},s\}]$ and $\bm{J}=\{0\}$, for some $i\in \{1,2,\cdots,\lceil4\cdot 3^{2/3}s\rceil\}$. The condition on the inclusion for the intervals $\bm{I}$ and $\bm{J}$ is ensured exactly in the same way as before, and the parabolic curvature is bounded above by $40s^2$. Hence, using Theorem \ref{t:unifmoddev} and a union bound, it follows that, for $s$ large enough and $n\geq 3n_0$,
\[\P(\mathsf{G}_2^c\cup \mathsf{G}_4^c)\leq \lceil 4\cdot 3^{2/3}s\rceil Ce^{-cs^3}\leq 6^{-1}\,.\]
Finally, for $\mathsf{G}_5$, observe that, since paths between two fixed endpoints constrained to stay in a box have smaller weight than does the polymer between these endpoints, we can again use Theorem \ref{t:unifmoddev}. For a given application of Theorem \ref{t:unifmoddev},  take $\bm{n}=n, \bm{t_1}=1/3,\bm{t_2}=2/3,
\bm{s}=40s^2,\bm{I}=[-s+(i-1)3^{-2/3},\max\{-s+i3^{-2/3},s\}]$ and $\bm{J}=[-s+(j-1)3^{-2/3},\max\{-s+j3^{-2/3},s\}]$ for $i\in \{1,2,\cdots,\lceil 2\cdot 3^{2/3}s\rceil\}$ and $j\in \{1,2,\cdots,\lceil 2\cdot 3^{2/3}s\rceil\}$. As before, the condition on inclusion for $\bm{I}$ and~$\bm{J}$ is satisfied, and the parabolic curvature is at most $3^{4/3}s^2$, which is less than $10s^2$. Thus, applying Theorem \ref{t:unifmoddev} and a union bound, we find that, for $n \geq 3n_0$ and $s$ large,
\[\P(\mathsf{G}_5^c)\leq \P\left(\sup \weight_{n;([-s,s],1/3)}^{([-s,s],2/3)}>50s^2\right)\leq \lceil 2\cdot 3^{2/3}s\rceil^2Ce^{-cs^3}\leq 6^{-1}\,.\]
Thus we have $\P(\mathsf{F}) \geq 1/2$ by a union bound.
In Case High, we also have 
$$
 \P \big( \high(\west,300s^2) \big) \geq \exp \big\{ -2C s^3 \big\} \, ,
$$
because $\west$ is a translate of $\midd$. Since the interior of $\west$  is disjoint from the regions that dictate the occurrence of $\mathsf{F}$, we see that 
\begin{equation}\label{e.fbound}
 \P \Big( \high(\west,300s^2) \cap \mathsf{F} \Big) \geq  2^{-1} \exp \big\{ -2C s^3 \big\} \, .
\end{equation}
When  $\high(\west,300s^2) \cap \mathsf{F}$ occurs, a high weight path connecting $(0,0)$ to $(0,1)$ may be formed by running it through $\west$. Indeed, and as Figure~\ref{f.casehigh} depicts, let $\rho$ denote a polymer running across, and contained in, $\west$,
whose weight is at least $300s^2$. If $x,y \in [-3s,-s]$  are such that $(x,1/3)$ and $(y,2/3)$ are $\rho$'s endpoints, then the path $\rho_{n;(0,0)}^{\leftarrow;(x,1/3)} \circ \rho \circ \rho_{n;(y,2/3)}^{\leftarrow,(0,1)}$
connects $(0,0)$ to $(0,1)$ and has weight at least $-50s^2 + 300s^2 - 50s^2$, in view of the first two conditions that specify $\mathsf{F}$.

On the other hand, the final three conditions specifying $\mathsf{F}$ ensure that, when this event occurs, any path from $(0,0)$ to $(0,1)$
whose $x$-coordinate never exceeds $s$ in absolute value has weight at most $50s^2 + 50s^2 + 50s^2$; indeed, the weight of any such path may be represented as a sum of the weights of the three subpaths formed by cutting the path at heights one-third and two-thirds. 

We thus find that, on $\high(\west,300s^2) \cap \mathsf{F}$, any path from $(0,0)$ to $(0,1)$ that remains in $\strip$ has weight at most $150s^2$; at the same time, a path of weight at least $200s^2$ connects these two points. Thus, we see that any polymer from $(0,0)$ to $(0,1)$ has maximum transversal fluctuation at least $s$ in this event. By~(\ref{e.fbound}), we find that
\begin{equation}\label{e.conchigh}
\P\left( \min \left\{ \mathrm{TF}(\rho) : \rho\in \Phi_{n;(0,0)}^{(0,1)} \right\}  \geq  s  \right) \,  \geq \,  2^{-1} \exp \big\{ -2C s^3 \big\} \, .
\end{equation}

Suppose now instead that Case Low holds.
We will argue that
\begin{equation}\label{e.tbargued}
 \P \Big( \weight_{n;(0,0)}^{(0,1)} \geq 1000 s^2 \, , \, \neg \,  \high\big( [-s,s] \times [0,1] , 900 s^2 \big)   \Big) \geq  2^{-1} \exp \big\{ - C s^3 \big\} \, ,
\end{equation}
where $\neg\, A$ denotes the complement of the event $A$.
Before we do so, we show that the event on this left-hand side entails that any polymer from $(0,0)$ to $(0,1)$ must leave the strip $[-s,s] \times [0,1]$; thus, (\ref{e.conchigh})
holds in Case Low, even when the factor of $2$ is omitted from the right-hand exponential.
When the last left-hand event occurs, any path from $(0,0)$ to $(0,1)$ that remains in the strip has weight at most $900s^2$. At the same time, the weight of any polymer from $(0,0)$ to $(0,1)$
is at least $1000s^2$. It is thus impossible for any polymer to remain in the strip.

To derive~(\ref{e.tbargued}), note that, because $\north$ and $\south$ are translates of $\midd$, Case Low entails that
$$
 \P \Big( \high(\south,300s^2) \cup  \high(\midd,300s^2) \cup  \high(\north,300s^2)  \Big) < 3 \exp \big\{ -2C s^3 \big\} \, .
$$
The bound~(\ref{e.wthou}) then yields~(\ref{e.tbargued}), since $3 \exp \big\{ -2 C s^3 \big\}  \leq 2^{-1} \exp \big\{ - C s^3 \big\}$ for all $s$ large enough.

The bound~(\ref{e.conchigh}) has been derived in both of the cases, so that proof of Proposition~\ref{p:lowbndTF} is complete. \qed

\bibliography{HolderLPP}
\bibliographystyle{plain}
\end{document}